\pdfoutput=1
\RequirePackage{ifpdf}
\ifpdf 
\documentclass[pdftex]{sigma}
\else
\documentclass{sigma}
\fi

\begin{document}

\allowdisplaybreaks

\renewcommand{\PaperNumber}{044}

\FirstPageHeading

\ShortArticleName{Vector Polynomials and a~Matrix Weight Associated to Dihedral Groups}

\ArticleName{Vector Polynomials and a~Matrix Weight\\
Associated to Dihedral Groups}

\Author{Charles F.~DUNKL}

\AuthorNameForHeading{C.F.~Dunkl}

\Address{Department of Mathematics, University of Virginia,\\
PO Box 400137, Charlottesville VA 22904-4137, USA}
\Email{\href{mailto:cfd5z@virginia.edu}{cfd5z@virginia.edu}}
\URLaddress{\url{http://people.virginia.edu/~cfd5z/}}

\ArticleDates{Received January 22, 2014, in f\/inal form April 10, 2014; Published online April 15, 2014}

\Abstract{The space of polynomials in two real variables with values in a~2-dimensional irreducible module of a~dihedral
group is studied as a~standard module for Dunkl operators.
The one-parameter case is considered (omitting the two-parameter case for even dihedral groups).
The matrix weight function for the Gaussian form is found explicitly by solving a~boundary value problem, and then
computing the normalizing constant.
An orthogonal basis for the homogeneous harmonic polynomials is constructed.
The coef\/f\/icients of these polynomials are found to be balanced terminating ${}_4F_3$-series.}


\Classification{33C52; 20F55; 33C45}

\section{Introduction}

For each irreducible two-dimensional representation of a~dihedral ref\/lection group there is a~mo\-dule of the algebra of
operators on polynomials generated by multiplication and the Dunkl ope\-rators.
This algebra is called the rational Cherednik algebra of the group.
The space of vector-valued polynomials is equipped with a~bilinear form which depends on one parameter and is invariant
under the group action.
For a~certain interval of parameter values the form can be represented as an integral with respect to
a~positive-def\/inite matrix weight function times the Gaussian measure.
In a~previous paper this structure was analyzed for the group of type~$B_{2}$ where there are two free parameters.
This paper concerns the Coxeter group of type $I_{2} ( m ) $, the full symmetry group of the regular $m$-gon,
of order~$2m$.
Even values for $m$ would allow two parameters but only the one parameter case is considered here.

The orthogonality properties of the Gaussian form are best analyzed by means of harmonic homogeneous polynomials.
These are studied in Section~\ref{section2}.
An inductive approach is used to produce the def\/initions.
Closed forms are obtained for these polynomials in the complex coordinate system.
The formulae are explicit but do require double sums, that is, the coef\/f\/icients are given as terminating balanced
${}_{4}F_{3}$-series.
Section~\ref{section3} contains the construction of an orthogonal basis of harmonic polynomials and the structure constants with
respect to the Gaussian form.
The results in Sections~\ref{section2} and~\ref{section3} are of algebraic f\/lavor and hold for any value of~$\kappa$.
The sequel is of analytic nature and the results hold for restricted values of~$\kappa$.
Section~\ref{section4} sets up the dif\/ferential system for the weight function and then constructs the solution up to a~normalizing
constant.
Section~\ref{section5} deals with the mechanics of integrating harmonic polynomials with respect to the Gaussian weight.
In Section~\ref{section6} the normalizing constant is found by direct integration.
The survey~\cite{Dunkl2012} covering general information about Dunkl operators and related topics is accessible online.

Fix $m\geq3$ and set $\omega:=\exp\frac{2\pi\mathrm{i}}{m}$.
The group $I_{2}(m) $ (henceforth denoted by $W$) contains $m$ ref\/lections and $m-1$ rotations.
Although the group is real it is often useful to employ complex coordinates for $x= ( x_{1},x_{2} )
\in\mathbb{R}^{2}$, that is, $z=x_{1}+\mathrm{i}x_{2}$, $\overline{z}=x_{1}-\mathrm{i}x_{2}$.
Then the ref\/lections are expressed as $z\sigma_{j}:=\overline{z}\omega^{j}$ ($0\leq j<m$) and the rotations are
$z\varrho_{j}:=z\omega^{j}$ $( 1\leq j<m ) $ The (pairwise nonequivalent) $2$-dimensional irreducible
representations are given by
\begin{gather*}
\tau_{\ell} ( \sigma_{j} ) =\left[
\begin{matrix}
0 & \omega^{-j\ell}
\\
\omega^{j\ell} & 0
\end{matrix}
\right],
\qquad
1\leq\ell\leq\left\lfloor \frac{m-1}{2}\right\rfloor,
\qquad
0\leq j<m.
\end{gather*}
Henceforth f\/ix $\ell$.
We consider the vector space $\mathcal{P}_{m,\ell}$ consisting of polynomials
\begin{gather*}
f ( z,\overline{z},t,\overline{t} ) =f_{1} ( z,\overline {z} ) t+f_{2} ( z,\overline{z} )
\overline{t}
\end{gather*}
with the group action
\begin{gather}
\sigma_{j}f ( z,\overline{z},t,\overline{t} ) =f_{2}\big( \overline{z}\omega^{j},z\omega^{-j}\big)
\omega^{-\ell j}t+f_{1}\big( \overline{z}\omega^{j},z\omega^{-j}\big) \omega^{\ell j}\overline{t}.
\label{grpact}
\end{gather}
Fix a~parameter $\kappa$; the Dunkl operators are def\/ined by
\begin{gather}
\mathcal{D}f ( z,\overline{z},t,\overline{t} )=\frac{\partial }{\partial z}f (
z,\overline{z},t,\overline{t} ) +\kappa\sum _{j=0}^{m-1}\frac{f ( z,\overline{z},\omega^{\ell
j}\overline{t},\omega^{-\ell j}t ) -f ( \overline{z}\omega^{j},z\omega ^{-j},\omega^{\ell
j}\overline{t},\omega^{-\ell j}t ) }{z-\overline {z}\omega^{j}},
\label{DDdef}
\\
\overline{\mathcal{D}}f ( z,\overline{z},t,\overline{t} )=\frac{\partial}{\partial\overline{z}}f (
z,\overline{z},t,\overline {t} ) -\kappa\sum\limits_{j=0}^{m-1}\frac{f ( z,\overline{z},\omega^{\ell
j}\overline{t},\omega^{-\ell j}t ) -f ( \overline{z}\omega ^{j},z\omega^{-j},\omega^{\ell
j}\overline{t},\omega^{-\ell j}t ) }{z-\overline{z}\omega^{j}}\omega^{j}.
\nonumber
\end{gather}
In these coordinates the Laplacian $\Delta_{\kappa}=4\mathcal{D} \overline{\mathcal{D}}$.
The group covariant property $\overline{\mathcal{D} }=\sigma_{0}\mathcal{D\sigma}_{0}$ is convenient for simplifying
some proofs.

The \textit{rational Cherednik algebra} associated with the data $ ( I_{2}(m),\kappa ) $ is the
abstract algebra formed from the product $\mathbb{C} [ z,\overline{z} ] \otimes\mathbb{C} [
\mathcal{D},\overline{\mathcal{D}} ] \otimes\mathbb{C}I_{2}(m) $ (polynomials in the two sets of
variables and the group algebra of $I_{2}(m) $) with relations like $\overline{\mathcal{D}}
=\sigma_{0}\mathcal{D\sigma}_{0}$ and $\overline{z}=\sigma_{0}z\sigma_{0}$, and commutations described in
Proposition~\ref{dzf}.
The algebra acts on~$\mathcal{P}_{m,\ell}$ by multiplication, and equations~\eqref{DDdef},~\eqref{grpact} for the
respective three components.
Thus there is a~representation of the algebra as operators on~$\mathcal{P}_{m,\ell}$, and so~$\mathcal{P}_{m,\ell}$ is
called the \textit{standard module} with lowest degree component of isotype $\tau_{\ell}$ (that is, the representation
$\tau_{\ell}$ of $I_{2}(m) $).

The bilinear form on $\mathcal{P}_{m,\ell}$ is def\/ined in the real coordinate system
($t=s_{1}+\mathrm{i}s_{2}$, $\mathcal{D}_{1}=\mathcal{D+}\overline
{\mathcal{D}},\allowbreak\mathcal{D}_{2}=\mathrm{i} ( \mathcal{D-} \overline{\mathcal{D}} ) $) by
\begin{gather*}
 \langle s_{j},s_{k} \rangle=\delta_{jk},
\qquad
 \langle x_{j}p ( x,s ),q ( x,s ) \rangle= \langle p (
x,s ),\mathcal{D}_{j}q ( x,s )  \rangle,
\qquad
 \langle s_{j},q ( x,s )  \rangle= \langle s_{j},q ( 0,s )  \rangle.
\end{gather*}
This def\/inition is equivalent to
\begin{gather*}
\langle p_{1}(x) s_{1}+p_{2}(x) s_{2},q(x,s) \rangle =\langle
s_{1},( p_{1}( \mathcal{D}_{1},\mathcal{D}_{2}) q) ( 0,s) \rangle +\langle
s_{2},( p_{2}( \mathcal{D}_{1},\mathcal{D}_{2}) q) (0,s) \rangle.
\end{gather*}
The form is real-bilinear and symmetric (see~\cite{Dunkl2013}).
To transform the def\/inition to the complex coordinate system we impose the conditions
\begin{gather*}
 \langle cp,q \rangle =\overline{c} \langle p,q \rangle,
\qquad
 \langle p,cq \rangle =c \langle p,q \rangle,
\qquad
c\in \mathbb{C},
\end{gather*}
thus $ \langle t,t \rangle = \langle s_{1}+\mathrm{i}s_{2},s_{1}+\mathrm{i}s_{2} \rangle
=2$, $ \langle t,\overline{t}  \rangle =0$, $ \langle \overline{t},\overline{t} \rangle =2$.
Now suppose $p$ is given in complex form $p ( z,\overline{z} ) =\sum\limits_{j,k}a_{jk}z^{j}\overline{z}^{k}$
then set $p^{\ast} ( z,\overline{z} ) =\sum\limits_{j,k}\overline{a_{jk}}z^{k}\overline{z}^{j}$ and replace
$x_{j}$ by $\mathcal{D}_{j}$ and conjugate.
This leads to $p^{\ast } ( \mathcal{D}_{1}+\mathrm{i}\mathcal{D}_{2},\mathcal{D}_{1}
-\mathrm{i}\mathcal{D}_{2} ) =p^{\ast} ( 2\overline{\mathcal{D} },2\mathcal{D} ) $ (in more detail, if
$p=\sum\limits_{j,k}a_{jk}z^{j} \overline{z}^{k}$ then $p^{\ast} ( 2\overline{\mathcal{D}},2\mathcal{D}  )
=\sum\limits_{j,k}\overline{a_{jk}}2^{j+k}\mathcal{D}^{j}\overline {\mathcal{D}}^{k}$).

Then the form is given by
\begin{gather}
 \langle p_{1} ( z,\overline{z} ) t+p_{2} ( z,\overline {z} ) \overline{t},q (
z,\overline{z},t,\overline{t} )  \rangle = \langle t,p_{1}^{\ast} ( 2\overline{\mathcal{D}
},2\mathcal{D} ) q \rangle |_{z=0}+ \langle \overline{t},p_{2}^{\ast} (
2\overline{\mathcal{D}},2\mathcal{D} ) q \rangle |_{z=0}.
\label{pzform}
\end{gather}

The form is def\/ined for all $\kappa$.
It is positive-def\/inite provided that $-\frac{\ell}{m}<\kappa<\frac{\ell}{m}$, as will be shown.
This is related to the results of Etingof and Stoica~\cite[Proposition~4.3]{Etingof2009} on unitary representations.
The (abstract) Gaussian form is def\/ined for $f,g\in \mathcal{P}_{m,\ell}$ by
\begin{gather*}
 \langle f,g \rangle _{G}= \big\langle e^{\Delta_{\kappa} /2}f,e^{\Delta_{\kappa}/2}g \big\rangle
= \big\langle e^{2\mathcal{D} \overline{\mathcal{D}}}f,e^{2\mathcal{D}\overline{\mathcal{D}}}g \big\rangle.
\end{gather*}
This has the important property that $ \langle pf,g \rangle _{G}= \langle f,p^{\ast}g \rangle _{G}$
for any scalar polynomial $p ( z,\overline{z} )$.
This motivates the attempt to express the pairing in integral form, specif\/ically
\begin{gather}
 \langle f,g \rangle _{G}=\int_{\mathbb{C}}gKf^{\ast}e^{- \vert z \vert ^{2}/2}dm_{2} ( z ),
\label{defK}
\end{gather}
where $f,g\in\mathcal{P}_{m,\ell}$ are interpreted as vectors $f_{1} ( z,\overline{z} ) t+f_{2} (
z,\overline{z} ) \overline {t}\mapsto ( f_{1},f_{2} ) $ and $K$ is an integrable Hermitian $2\times2$
matrix function.
We preview the formula for $K$ on the fundamental chamber $ \big\{
z=re^{\mathrm{i}\theta}:r>0,\, 0<\theta<\frac{\pi}{m}\big\} $, in the real coordinate system, in terms of auxiliary
parameters and functions:

Set $v:=\sin^{2}\frac{m\theta}{2}$, $\delta:=\frac{1}{2}-\frac{\ell}{m}$,
\begin{gather*}
f_{1} ( \kappa,\delta;v ) :=v^{\kappa/2} ( 1-v ) ^{-\kappa/2}\;{}_{2}F_{1}\left(
\genfrac{}{}{0pt}{}{\delta,-\delta}{\frac{1}{2}+\kappa} ;v\right),
\\
f_{2} ( \kappa,\delta;v ) :=\frac{\delta}{\frac{1}{2}+\kappa }v^{ ( \kappa+1 ) /2} ( 1-v )
^{ ( 1-\kappa ) /2}\;{}_{2}F_{1}\left(   \genfrac{}{}{0pt}{}{1+\delta,1-\delta}{\frac{3}{2}+\kappa} ;v\right),
\\
H ( \kappa,\delta ) :=\frac{\Gamma\left( \frac{1}{2} +\kappa\right) ^{2}}{\Gamma\left(
\frac{1}{2}+\kappa+\delta\right) \Gamma\left( \frac{1}{2}+\kappa-\delta\right) },
\end{gather*}
and
\begin{gather*}
L ( \theta ) =\left[
\begin{matrix}
f_{1} ( \kappa,\delta;v ) & f_{2} ( \kappa,\delta;v )
\\
-f_{2} ( -\kappa,\delta;v ) & f_{1} ( -\kappa,\delta;v )
\end{matrix}
\right] \left[
\begin{matrix}
\cos\delta m\theta & -\sin\delta m\theta
\\
\sin\delta m\theta & \cos\delta m\theta
\end{matrix}
\right],
\end{gather*}
then
\begin{gather*}
K ( r,\theta ) =\frac{\cos\pi\kappa}{2\pi\cos\pi\delta}L ( \theta ) ^{T}\left[
\begin{matrix}
H ( -\kappa,\delta ) & 0
\\
0 & H ( \kappa,\delta )
\end{matrix}
\right] L ( \theta ).
\end{gather*}
The equation $K ( zw ) =\tau_{\ell} ( w ) ^{-1}K ( z ) \tau_{\ell} ( w ) $ def\/ines
$K$ on the other chambers  $\big\{ z=re^{\mathrm{i}\theta}:r>0$, $( j-1 ) \frac{\pi}{m}
<\theta<j\frac{\pi}{m}\big\}$, $2\leq j\leq2m$.

Usually we will write $\partial_{u}:=\frac{\partial}{\partial u}$ for a~variable $u$.

\section{Harmonic polynomials}  \label{section2}

In this section we construct a~basis for the (vector-valued) harmonic ($\Delta_{\kappa}f=0$) homogeneous polynomials.
These are important here because for generic $\kappa$ any polynomial in $\mathcal{P}_{m,\ell}$ has a~unique expansion as
a~sum of terms like $ ( z\overline{z} ) ^{n}f ( z,\overline{z},t,\overline{t} ) $ where~$f$ is
homogeneous and $\Delta_{\kappa}f=0$ (see~\cite[(5),~p.~4]{Dunkl2013}); furthermore $ \langle f,g \rangle
_{G}= \langle f,g \rangle $ for any harmonic polynomials~$f$,~$g$.
Complex coordinates are best suited for this study.
To each monomial in~$z$,~$\overline{z}$, $t$ or~$\overline{t}$ associate the degree and the $m$-parity: (with
$a,b\in\mathbb{N}_{0}$){\samepage
\begin{enumerate}\itemsep=0pt
\item[1)] the degree is given by $\deg ( z^{a}\overline{z}^{b}t ) =\deg ( z^{a}\overline{z}^{b}\overline{t} )
=a+b$,
\item[2)] the $m$-parities of $z^{a}\overline{z}^{b}t$, $z^{a}\overline{z} ^{b}\overline{t}$ are $a-b+\ell\operatorname{mod}m$,
$a-b-\ell \operatorname{mod}m$, respectively.
\end{enumerate}}

If each monomial in a~polynomial $p$ has the same degree and $m$-parity then say that~$p$ has this degree and $m$-parity
(this implies~$p$ is homogeneous in~$z$,~$\overline{z}$).
The following is a~product rule.

\begin{proposition}\label{dzf}
For $f\in\mathcal{P}_{m,\ell}$
\begin{gather*}
\mathcal{D}zf ( z,\overline{z},t,\overline{t} )= ( z\mathcal{D}+1 ) f (
z,\overline{z},t,\overline{t} ) +\kappa\sum\limits_{j=0}^{m-1}f\big(
\overline{z}\omega^{j},z\omega^{-j},\omega^{\ell j}\overline{t},\omega^{-\ell j}t\big),
\\
\overline{\mathcal{D}}zf ( z,\overline{z},t,\overline{t} )=z\overline{\mathcal{D}}f (
z,\overline{z},t,\overline{t} ) -\kappa\sum\limits_{j=0}^{m-1}f\big(
\overline{z}\omega^{j},z\omega^{-j},\omega^{\ell j}\overline{t},\omega^{-\ell j}t\big) \omega^{j}.
\end{gather*}
\end{proposition}

\begin{proof}
Abbreviate $g_{j}=f ( z,\overline{z},\omega^{\ell j}\overline{t},\omega^{-\ell j}t ) $ and $h_{j}=f (
\overline{z}\omega ^{j},z\omega^{-j},\omega^{\ell j}\overline{t},\omega^{-\ell j}t )$.
Then
\begin{gather*}
\mathcal{D}zf ( z,\overline{z},t,\overline{t} )=\left( z\frac{\partial}{\partial z}+1\right) f (
z,\overline{z},t,\overline {t} ) +\kappa\sum\limits_{j=0}^{m-1}\frac{zg_{j}-\overline{z}\omega^{j}h_{j}
}{z-\overline{z}\omega^{j}}
\\
\hphantom{\mathcal{D}zf ( z,\overline{z},t,\overline{t} )}{}
=\left( z\frac{\partial}{\partial z}+1\right) f ( z,\overline {z},t,\overline{t} )
+\kappa\sum\limits_{j=0}^{m-1}\left\{ \frac{zg_{j}
-zh_{j}}{z-\overline{z}\omega^{j}}+\frac{zh_{j}-\overline{z}\omega^{j}h_{j} }{z-\overline{z}\omega^{j}}\right\}
\\
\hphantom{\mathcal{D}zf ( z,\overline{z},t,\overline{t} )}{}
= ( z\mathcal{D}+1 ) f ( z,\overline{z},t,\overline {t} ) +\kappa\sum\limits_{j=0}^{m-1}h_{j}.
\end{gather*}
The second formula is proved similarly.
\end{proof}

With the aim of using Proposition~\ref{dzf} def\/ine
\begin{gather*}
T_{z}f ( z,\overline{z},t,\overline{t} ) :=\sum\limits_{j=0} ^{m-1}f\big(
\overline{z}\omega^{j},z\omega^{-j},\omega^{\ell j}\overline {t},\omega^{-\ell j}t\big),
\\
T_{\overline{z}}f ( z,\overline{z},t,\overline{t} ) :=\sum\limits_{j=0}^{m-1}f\big(
\overline{z}\omega^{j},z\omega^{-j},\omega^{\ell j}\overline{t},\omega^{-\ell j}t\big) \omega^{j}.
\end{gather*}

\begin{lemma}\label{Tzzb}
Suppose that $p\in\mathcal{P}_{m,\ell}$ has $m$-parity $r$ then $T_{z}p=0$ for $r\neq0\operatorname{mod}m$,
$T_{z}p=m\sigma_{0}p$ for $r=0\operatorname{mod}m$, $T_{\overline{z}}p=0$ for $r\neq m-1\operatorname{mod}m$ and
$T_{\overline{z}}p=m\sigma_{0}p$ if $r=m-1\operatorname{mod}m$.
\end{lemma}

\begin{proof}
Consider $T_{z}z^{a}\overline{z}^{b}t=z^{b}\overline{z}^{a}\overline{t} \sum\limits_{j=0}^{m-1}\omega^{ (
a-b+\ell ) j}$; the sum is zero if $a-b+\ell\neq0\operatorname{mod}m$, otherwise $T_{z}z^{a}\overline{z}
^{b}t=m\sigma_{0} ( z^{a}\overline{z}^{b}t ) $.
The same argument applies to $T_{z}z^{a}\overline{z}^{b}\overline{t}$ whose $m$-parity is $a-b-\ell$.
Similarly $T_{\overline{z}}z^{a}\overline{z}^{b}t=z^{b}
\overline{z}^{a}\overline{t}\sum\limits_{j=0}^{m-1}\omega^{ ( a-b+\ell+1 ) j}$.
\end{proof}

A polynomial $p$ is called harmonic if $\mathcal{D}\overline{\mathcal{D}}p=0$.
By using an inductive construction we do not need to explicitly determine the action of
$\mathcal{D}$, $\overline{\mathcal{D}}$ on general monomials.

\begin{proposition}\label{Dzp}
Suppose that $p\in\mathcal{P}_{m,\ell}$ has $m$-parity $r$ with $r\neq0\operatorname{mod}m$, and
$\overline{\mathcal{D}}p=0$ then $\mathcal{D}zp= ( z\mathcal{D}+1 ) p$ and $\overline{\mathcal{D} }zp=-\kappa
m\sigma_{0}p$ if $r=m-1\operatorname{mod}m$, otherwise $\overline{\mathcal{D}}zp=0$.
\end{proposition}

\begin{proof}
By Proposition~\ref{dzf} $\mathcal{D}zp= ( z\mathcal{D}+1 ) p+\kappa T_{z}p$ and $T_{z}p=0$ by the lemma.
Also $\overline{\mathcal{D} }zp=z\overline{\mathcal{D}}p-\kappa T_{\overline{z}}p$.
\end{proof}

Notice that $p=t$ or $p=\overline{t}$ with degree $0$ and $m$-parities $\ell$ and $m-\ell$ respectively, satisfy the
hypotheses of the following.

\begin{theorem}\label{p2zp}
Suppose that $p\in\mathcal{P}_{m,\ell}$ has $m$-parity $r$ with $1\leq r\leq m-1$, $\deg p=n$, $z\mathcal{D}p=np$ and
$\overline{\mathcal{D} }p=0$ then $\mathcal{D}z^{s}p= ( n+s ) z^{s-1}p$ for $0\leq s\leq m-r$,
$\overline{\mathcal{D}}z^{s}p=0$ for $0\leq s\leq m-1-r$ and
$\overline{\mathcal{D}}z^{m-r}p=-m\kappa\overline{z}^{m-1-r}\sigma_{0}p$.
Furthermore $q:=z^{m-r}p+\frac{m\kappa}{m-r+n}\overline{z}^{m-r}\sigma_{0}p$ has $m$-parity $0$ and satisfies
\begin{gather*}
\mathcal{D}q= ( m-r+n ) \left\{ 1-\left( \frac{m\kappa }{m-r+n}\right) ^{2}\right\}
z^{m-r-1}p,\qquad \overline{\mathcal{D}}q=0,
\\
\mathcal{D}zq= ( m+n-r+1 ) q,\qquad \overline{\mathcal{D}}zq=0.
\end{gather*}
\end{theorem}

\begin{proof}
By Lemma~\ref{Tzzb} $T_{z}z^{s}p=0$ for $0\leq s<m-r$.
Arguing by induction suppose $\mathcal{D}z^{s}p= ( n+s ) z^{s-1}p$ for some $s$ satisfying $1\leq s<m-r$ then
by Proposition~\ref{Dzp} $\mathcal{D} z^{s+1}p= ( z\mathcal{D}+1 ) z^{s}p= (  ( n+s ) +1 )
z^{s}p$.
Similarly $T_{\overline{z}}z^{s}p=0$ for $0\leq s\leq m-2-r$ and $\overline{\mathcal{D}}z^{s}p=0$ for $1\leq s\leq
m-1-r$.
Also $\overline{\mathcal{D}}z^{m-r}p=-m\kappa\sigma_{0} ( z^{m-r-1}p ) $.
Thus
\begin{gather*}
\overline{\mathcal{D}}q=\overline{\mathcal{D}}z^{m-r}p+\frac{m\kappa
}{m\!-\!r\!+\!n}\overline{\mathcal{D}}\mathcal{\sigma}_{0}\big( z^{m-r}p\big)
=-m\kappa\sigma_{0}\big( z^{m-r-1}p\big) +\frac{m\kappa}{m\!-\!r\!+\!n} \sigma_{0}\mathcal{D}z^{m-r}p=0,
\end{gather*}
and
\begin{gather*}
\mathcal{D}q=Dz^{m-r}p+\frac{m\kappa}{m-r+n}\mathcal{D\sigma}_{0}\big( z^{m-r}p\big)
\\
\hphantom{\mathcal{D}q}{}
=( m-r+n) z^{m-r-1}p+\frac{m\kappa}{m-r+n}\sigma_{0} \overline{\mathcal{D}}z^{m-r}p
\\
\hphantom{\mathcal{D}q}{}
= ( m-r+n ) z^{m-r-1}p+\frac{m\kappa}{m-r+n}\sigma_{0}\big( -m\kappa\sigma_{0}z^{m-r-1}p\big)
\\
\hphantom{\mathcal{D}q}{}
= ( m-r+n ) \left\{ 1-\left( \frac{m\kappa}{m-r+n}\right) ^{2}\right\} z^{m-r-1}p.
\end{gather*}
The $m$-parity of $q$ is zero and $\overline{\mathcal{D}}q=0$ thus $\overline{\mathcal{D}}zq=0$.
By Proposition~\ref{dzf}
\begin{gather*}
\mathcal{D}zq=z\mathcal{D}q+q+m\kappa\sigma_{0}q= ( m-r+n ) \left\{ 1-\left( \frac{m\kappa}{m-r+n}\right)
^{2}\right\} z^{m-r}p
\\
\hphantom{\mathcal{D}zq=}
{}+z^{m-r}p+\frac{m\kappa}{m-r+n}\overline{z}^{m-r}\sigma_{0}p+m\kappa \left\{
\overline{z}^{m-r}\sigma_{0}p+\frac{m\kappa}{m-r+n}z^{m-r}p\right\}
\\
\hphantom{\mathcal{D}zq }{}
= ( m-r+n+1 ) \left\{ z^{m-r}p+\frac{m\kappa}{m-r+n} \overline{z}^{m-r}\sigma_{0}p\right\}.
\end{gather*}
This completes the proof.
\end{proof}

\begin{corollary}
With the same hypotheses $z^{s}p$, $\overline{z}^{s}\sigma_{0}p$ for $0\leq s\leq m-r$, and $zq$, $\sigma_{0} ( zq )
$ are harmonic.
\end{corollary}

\begin{proof}
Each of these polynomials $f$ satisfy $\mathcal{D}f=0$ or $\overline {\mathcal{D}}f=0$ (recall
$\mathcal{D\sigma}_{0}=\sigma_{0}\overline {\mathcal{D}}$) with the exception of $z^{m-r}p$ and $\sigma_{0} (
z^{m-r}p ) $.
But $\overline{\mathcal{D}}\mathcal{D}z^{m-r} p=  ( m-r+n ) \overline{\mathcal{D}}z^{m-r-1}p=0$.
\end{proof}

Observe that $\mathcal{D}$, $\overline{\mathcal{D}}$ change the $m$-parity by $-1,+1$ respectively.
The \textit{leading term} of a~homogeneous polynomial
$\sum\limits_{j=0}^{n}a_{j}z^{n-j}\overline{z}^{j}t+\sum\limits_{j=0}^{n}b_{j} z^{n-j}\overline{z}^{j}\overline{t}$ is
def\/ined to be $ ( a_{0} z^{n}+a_{n}\overline{z}^{n} ) t+ ( b_{0}z^{n}+b_{n}\overline{z} ^{n} )
\overline{t}$.

There are $4$ linearly independent harmonic polynomials of each degree ${\geq}1$.
Recall  $\sigma_{0}p ( z,\overline{z},t,\overline{t} ) $ $=p ( \overline{z},z,\overline{t},t )
$.
By use of Theorem~\ref{p2zp} we f\/ind there are two sequences of harmonic homogeneous polynomials $ \big\{ p_{n}^{ (
1 ) },\sigma_{0}p_{n}^{ ( 1 ) }:n\geq 1 \big\} $ and $ \big\{ p_{n}^{ ( 2 )
},\sigma_{0}p_{n}^{ ( 2 ) }:n\geq1 \big\} $ with the properties:
\begin{enumerate}\itemsep=0pt
\item[1)]
the leading terms of $p_{n}^{(1) }$, $\sigma_{0} p_{n}^{(1) }$, $p_{n}^{(2)
}$, $\sigma_{0}p_{n}^{(2) }$ are $z^{n}t$, $\overline{z}^{n}\overline{t}$, $z^{n}\overline {t}$, $\overline{z}^{n}t$
respectively,

\item[2)]
$p_{1}^{(1) }=zt$, $\sigma_{0}p_{1}^{(1) }=\overline{z}\overline{t}$,

\item[3)]
if $n+\ell\neq0\operatorname{mod}m$ then $p_{n+1}^{(1) }=zp_{n}^{(1) }$,

\item[4)]
if $n+\ell=0\operatorname{mod}m$ then $p_{n+1}^{(1) }=z\big( p_{n}^{(1)
}+\frac{m\kappa}{n}\sigma_{0} p_{n}^{(1) }\big) $;

\item[5)]
$p_{1}^{(2) }=z\overline{t}$, $\sigma_{0}p_{1}^{(2) }=\overline{z}t$,

\item[6)]
if $n-\ell\neq0\operatorname{mod}m$ then $p_{n+1}^{(2) }=zp_{n}^{(2) }$,

\item[7)]
if $n-\ell=0\operatorname{mod}m$ then $p_{n+1}^{(2) }=z\big( p_{n}^{(2)
}+\frac{m\kappa}{n}\sigma_{0} p_{n}^{(2) }\big) $.
\end{enumerate}

The consequence of these relations is that there exist polynomials $P_{n}^{(1) }$ and $P_{n}^{(
2) }$ of degrees $m( n+1) -\ell+1$ and $nm+\ell+1$ respectively ($n\geq0$) such that
\begin{enumerate}\itemsep=0pt
\item[1)]
$p_{r}^{(1) }=z^{r}t$ for $1\leq r\leq m-\ell$ and $p_{s}^{(1) }=z^{r}P_{n}^{(1) }$
for $s=m(n+1) -\ell+1+r$ and $0\leq r<m$;

\item[2)]
$p_{r}^{(2) }=z^{r}\overline{t}$ for $1\leq r\leq\ell$ and $p_{s}^{(2) }=z^{r}P_{n}^{(2)}$ for $s=nm+\ell+1+r$ and $0\leq r<m$.
\end{enumerate}

The formulae imply recurrences:
\begin{enumerate}\itemsep=0pt
\item[1)]
$P_{0}^{(1) }=z^{m-\ell+1}t+\frac{m\kappa}{m-\ell }z\overline{z}^{m-\ell}\overline{t}$ and $P_{n+1}^{(1) } =z^{m}P_{n}^{(1) }+\frac{m\kappa}{(n+1) m
+(m-\ell)
}z\overline{z}^{m-1}\sigma_{0}P_{n}^{(1) }$;

\item[2)]
$P_{0}^{(2) }=z^{\ell+1}\overline{t}+\frac{m\kappa}{\ell }z\overline{z}^{\ell}t$ and $P_{n+1}^{(2) }=z^{m}P_{n}^{(2) }+\frac{m\kappa}{(n+1) m+\ell}z\overline{z}
^{m-1}\sigma_{0}P_{n}^{(2) }$.
\end{enumerate}

There is a~common thread in these recurrences.
Let $\lambda$ be a~parameter with $\lambda>0$ and def\/ine polynomials $Q_{n}^{(1) } (
\kappa,\lambda;w,\overline{w} ) $, $Q_{n}^{(2) } ( \kappa,\lambda;w,\overline{w} ) $ by
\begin{gather*}
Q_{0}^{(1) } ( \kappa,\lambda;w,\overline{w} )=1,
\qquad
Q_{0}^{(2) } (\kappa,\lambda;w,\overline{w} ) =\frac{\kappa}{\lambda},
\\
Q_{n+1}^{(1) } ( \kappa,\lambda;w,\overline{w} )=wQ_{n}^{(1) } (
\kappa,\lambda;w,\overline{w} ) +\frac{\kappa}{\lambda+n+1}\overline{w}Q_{n}^{(2) } (
\kappa,\lambda;\overline{w},w ),
\\
Q_{n+1}^{(2) } ( \kappa,\lambda;w,\overline{w} )=\frac{\kappa}{\lambda+n+1}\overline{w}Q_{n}^{(1) } ( \kappa,\lambda;\overline{w},w )
+wQ_{n}^{(2) } ( \kappa,\lambda;w,\overline{w} ).
\end{gather*}

There is a~closed form point evaluation:

\begin{proposition}
For $\lambda>0$ and $n\geq0$
\begin{gather*}
Q_{n}^{(1) } ( \kappa,\lambda;1,1 )=\frac{ ( \lambda+\kappa ) _{n+1}+ (
\lambda-\kappa ) _{n+1}}{2 ( \lambda ) _{n+1}},
\\
Q_{n}^{(2) } ( \kappa,\lambda;1,1 )=\frac{ ( \lambda+\kappa ) _{n+1}- (
\lambda-\kappa ) _{n+1}}{2 ( \lambda ) _{n+1}}.
\end{gather*}
\end{proposition}

\begin{proof}
The statement is obviously true for $n=0$.
Arguing by induction suppose the statement is valid for some $n$ then
\begin{gather*}
Q_{n+1}^{(1) } ( \kappa,\lambda;1,1 )=\frac{ ( \lambda+\kappa ) _{n+1}+ (
\lambda-\kappa ) _{n+1}}{2(\lambda) _{n+1}}+\frac{\kappa}{\lambda+n+1}\left\{ \frac{(
\lambda+\kappa) _{n+1}-(\lambda-\kappa) _{n+1}}{2(\lambda) _{n+1}}\right\}
\\
\hphantom{Q_{n+1}^{(1) } ( \kappa,\lambda;1,1 )}{}
=\frac{(\lambda+\kappa) _{n+1}}{2(\lambda) _{n+1}}\left\{ 1+\frac{\kappa}{\lambda+n+1}\right\}
+\frac{(\lambda-\kappa) _{n+1}}{2(\lambda) _{n+1}}\left\{ 1-\frac{\kappa}{\lambda+n+1}\right\}
\\
\hphantom{Q_{n+1}^{(1) } ( \kappa,\lambda;1,1 )}{}
=\frac{(\lambda+\kappa) _{n+2}}{2(\lambda) _{n+2}}+\frac{(\lambda-\kappa)
_{n+2}}{2(\lambda) _{n+2}},
\end{gather*}
and
\begin{gather*}
Q_{n+1}^{(2) }( \kappa,\lambda;1,1)=\frac{\kappa}{\lambda+n+1}\left\{ \frac{(
\lambda+\kappa) _{n+1}+(\lambda-\kappa) _{n+1}}{2(\lambda) _{n+1} }\right\} +\frac{(
\lambda+\kappa) _{n+1}-( \lambda -\kappa) _{n+1}}{2(\lambda) _{n+1}}
\\
\hphantom{Q_{n+1}^{(2) }( \kappa,\lambda;1,1)}{}
=\frac{(\lambda+\kappa) _{n+1}}{2(\lambda) _{n+1}}\left\{ 1+\frac{\kappa}{\lambda+n+1}\right\}
-\frac{(\lambda-\kappa) _{n+1}}{2(\lambda) _{n+1}}\left\{ 1-\frac{\kappa}{\lambda+n+1}\right\}.
\end{gather*}
This completes the proof.
\end{proof}

\begin{proposition}
\label{PzQ}
For $n\geq0$
\begin{gather*}
P_{n}^{(1) }=z^{m-\ell+1}Q_{n}^{(1) }\left(
\kappa,\frac{m-\ell}{m};z^{m},\overline{z}^{m}\right) t+z\overline{z} ^{m-\ell}Q_{n}^{(2) }\left(
\kappa,\frac{m-\ell}{m};z^{m},\overline{z}^{m}\right) \overline{t},
\\
P_{n}^{(2) }=z^{\ell+1}Q_{n}^{(1) }\left( \kappa,\frac{\ell}{m};z^{m},\overline{z}^{m}\right)
\overline{t} +z\overline{z}^{\ell}Q_{n}^{(2) }\left( \kappa,\frac{\ell} {m};z^{m},\overline{z}^{m}\right) t.
\end{gather*}

\end{proposition}

\begin{proof}
The statement is true for $n=0$.
Suppose that the f\/irst one holds for some $n$, then set $\lambda=\frac{m-\ell}{m}$ and
\begin{gather*}
P_{n+1}^{(1) }=z^{m}P_{n}^{(1) }+\frac {m\kappa}{(n+1) m+(m-\ell)
}z\overline{z} ^{m-1}\sigma_{0}P_{n}^{(1) }
\\
\hphantom{P_{n+1}^{(1) }}{}
=z^{m}\left( z^{m-\ell+1}Q_{n}^{(1) }\big( \kappa,\lambda;z^{m},\overline{z}^{m}\big)
t+z\overline{z}^{m-\ell}Q_{n}^{(2) }\big( \kappa,\kappa;z^{m},\overline{z}^{m}\big) \overline {t} \right)
\\
\hphantom{P_{n+1}^{(1) }=}
{}+\frac{\kappa}{n+1+\lambda}z\overline{z}^{m-1}\left( \overline{z} ^{m-\ell+1}Q_{n}^{(1) }\big(
\kappa,\lambda;\overline{z} ^{m},z^{m} \big) \overline{t}+z^{m-\ell}\overline{z}Q_{n}^{(2) }\big(
\kappa,\lambda;\overline{z}^{m},z^{m}\big) t\right)
\\
\hphantom{P_{n+1}^{(1) }}{}
=z^{m-\ell+1}t\left( z^{m}Q_{n}^{(1) }\big( \kappa,\lambda;z^{m},\overline{z}^{m}\big)
+\frac{\kappa}{n+1+\lambda} \overline{z}^{m}Q_{n}^{(2) }\big( \kappa,\lambda;\overline {z}^{m},z^{m}\big)
\right)
\\
\hphantom{P_{n+1}^{(1) }=}
{}+z\overline{z}^{m-\ell}\overline{t}\left( z^{m}Q_{n}^{(2) }\big(
\kappa,\kappa;z^{m},\overline{z}^{m}\big) +\frac{\kappa }{n+1+\lambda}\overline{z}^{m}Q_{n}^{(1) }\big(
\kappa,\lambda;\overline{z}^{m},z^{m} \big) \right)
\\
\hphantom{P_{n+1}^{(1) }}{}
=z^{m-\ell+1}Q_{n+1}^{(1) }\big( \kappa,\lambda;z^{m},\overline{z}^{m}\big)
t+z\overline{z}^{m-\ell}Q_{n+1}^{(2) }\big( \kappa,\kappa;z^{m},\overline{z}^{m}\big) \overline{t};
\end{gather*}
this proves the f\/irst part by induction.
The same argument applies with $\ell$ replaced by $m-\ell$ and $t$, $\overline{t}$ interchanged.
\end{proof}

\begin{proposition}
For $n\geq0$
\begin{gather*}
Q_{n}^{(1) } ( \kappa,\lambda;w,\overline{w} )=w^{n}+\sum\limits_{j=1}^{n}\frac{\kappa^{2} (
n-j+1 ) }{\lambda ( \lambda+n ) }\; {}_{4}F_{3}\left(
\genfrac{}{}{0pt}{}{1-j,j-n,1-\kappa,1+\kappa}{2,\lambda+1,-\lambda-n+1} ;1\right) w^{n-j}\overline{w}^{j},
\\
Q_{n}^{(2) } ( \kappa,\lambda;w,\overline{w} )=\sum\limits_{j=0}^{n}\frac{\kappa}{\lambda+j}\; {}_{4}F_{3}\left(
\genfrac{}{}{0pt}{}{-j,j-n,-\kappa,+\kappa}{1,\lambda,-\lambda-n} ;1\right) w^{n-j}\overline{w}^{j}.
\end{gather*}
\end{proposition}

\begin{proof}
Set $Q_{n}^{(1) } ( \kappa,\lambda;w,\overline{w} )
=\sum\limits_{j=0}^{n}a_{n,j}w^{n-j}\overline{w}^{j}$ and $Q_{n}^{(2) } (
\kappa,\lambda;w,\overline{w} ) =\sum\limits_{j=0} ^{n}b_{n,j}w^{n-j}\overline{w}^{j}$, then the recurrence is
equivalent to $a_{0,0}=1$, $b_{0,0}=\frac{\kappa}{\lambda}$ and
\begin{gather*}
a_{n+1,j}=a_{n,j}+\frac{\kappa}{\lambda+n+1}b_{n,n+1-j},
\qquad
b_{n+1,j}=b_{n,j}+\frac{\kappa}{\lambda+n+1}a_{n,n+1-j},
\end{gather*}
for $0\leq j\leq n+1$ (with $a_{n,n+1}=0=b_{n,n+1}$).
It is clear that $a_{n+1,0}=a_{n,0}=1$ and $b_{n+1,0}=b_{n,0}=\frac{\kappa}{\lambda}$.
Now let $1\leq j\leq n+1$. In the following sums the upper limit can be taken as $n+1$.
Proceeding by induction we verify that $a_{n+1,j}-a_{n,j}=\frac{\kappa }{\lambda+n+1}b_{n,n+1-j}$ using the formulae.
The left-hand side is
\begin{gather*}
\frac{\kappa^{2} ( n-j+2 ) }{\lambda ( \lambda+n+1 ) }\;{}_{4}F_{3}\left(
\genfrac{}{}{0pt}{}{1-j,j-n-1,1-\kappa,1+\kappa}{2,\lambda+1,-\lambda-n} ;1\right)
\\
\qquad \quad {}-\frac{\kappa^{2} ( n-j+1) }{\lambda ( \lambda+n ) }\;{} _{4}F_{3}\left(
\genfrac{}{}{0pt}{}{1-j,j-n,1-\kappa,1+\kappa}{2,\lambda+1,-\lambda-n+1} ;1\right)
\\
\qquad {}
=\frac{\kappa^{2}}{\lambda}\sum\limits_{i=0}^{n+1}\! \frac{ ( 1-j ) _{i} ( 1-\kappa ) _{i} (
1+\kappa ) _{i}}{(2) _{i} ( \lambda+1 ) _{i}i!}\! \left\{ \frac{ ( n-j+2 )  (
j-n-1 ) _{i}}{ ( \lambda+n+1 )  ( -\lambda-n ) _{i}}-\frac{ ( n-j+1 )  ( j-n )
_{i} }{ ( \lambda+n )  ( -\lambda-n+1 ) _{i}}\right\}
\\
\qquad{} =\frac{\kappa^{2}}{\lambda}\sum\limits_{i=0}^{n+1}\frac{(1-j) _{i}(1-\kappa) _{i} (
1+\kappa ) _{i}}{ ( i+1 ) ! ( \lambda+1 ) _{i}i!}\left\{ \frac{ ( j-n-2 ) _{i+1}}{ (
-\lambda-n-1 ) _{i+1}}-\frac{ ( j-n-1 ) _{i+1}}{ ( -\lambda-n ) _{i+1}}\right\}
\\
\qquad{}
=\frac{\kappa^{2}}{\lambda}\sum\limits_{i=0}^{n+1}\frac{(1-j) _{i}(1-\kappa) _{i} (
1+\kappa ) _{i} ( j-n-1 ) _{i}}{ ( i+1 ) ! ( \lambda+1 ) _{i} i!}\frac{ ( i+1 )
 ( \lambda+j-1 ) }{ ( -\lambda-n-1 ) _{i+2}}
\\
\qquad{}
=\frac{\kappa^{2} ( \lambda+j-1 ) }{\lambda ( \lambda+n+1 ) ( \lambda+n)
}\;{}_{4}F_{3}\left(   \genfrac{}{}{0pt}{}{1-j,j-n-1,1-\kappa,1+\kappa}{1,\lambda+1,-\lambda-n+1} ;1\right).
\end{gather*}
The~$_{4}F_{3}$ is balanced and we can apply the transformation (see~\cite[(16.4.14)]{Olver2010}) to show
\begin{gather*}
{}_{4}F_{3}\left(   \genfrac{}{}{0pt}{}{1-j,j-n-1,1-\kappa,1+\kappa}{1,\lambda+1,-\lambda-n+1};1\right)
\\
\qquad{}
=\frac{ ( \lambda+2-j+n ) _{j-1} ( -\lambda+2-j ) _{j-1}}{ ( \lambda+1 ) _{j-1} (
-\lambda-n+1 ) _{j-1} }\;{}_{4}F_{3}\left(   \genfrac{}{}{0pt}{}{1-j,j-n-1,-\kappa,\kappa}{1,\lambda,-\lambda-n};1\right)
\\
\qquad{}
=\frac{\lambda ( \lambda+n ) }{ ( \lambda+j-1 )  ( \lambda+n+1-j )
}\frac{\lambda+n+1-j}{\kappa}b_{n,n-j+1},
\end{gather*}
and this proves the f\/irst recurrence relation.

Next verify $b_{n+1,j}-b_{n,j}=\frac{\kappa}{\lambda+n+1}a_{n,n+1-j}$.
The left-hand side equals
\begin{gather*}
\frac{\kappa}{\lambda+j}\sum\limits_{i=0}^{n+1}\frac{(-j) _{i}(-\kappa) _{i} (
\kappa ) _{i}}{i!i!(\lambda) _{i}}\left\{ \frac{ ( j-n-1 ) _{i}}{ ( -\lambda-n-1 )
_{i}}-\frac{ ( j-n ) _{i}}{ ( -\lambda-n ) _{i}}\right\}
\\
\qquad{}
=\frac{\kappa}{\lambda+j}\sum\limits_{i=1}^{n+1}\frac{(-j) _{i}(-\kappa) _{i} (
\kappa ) _{i} ( j-n ) _{i-1}}{i!i!(\lambda) _{i} ( -\lambda-n ) _{i-1} }\left\{
\frac{j-n-1}{-\lambda-n-1}-\frac{j-n+i-1}{-\lambda-n+i-1}\right\}
\\
\qquad{}
=\frac{\kappa}{\lambda+j}\sum\limits_{i=1}^{n+1}\frac{(-j) _{i}(-\kappa) _{i} (
\kappa ) _{i} ( j-n ) _{i-1} ( \lambda+j ) i}{i!i!(\lambda) _{i} (
-\lambda-n-1 ) _{i+1}}
\\
\qquad{}
=\left( \frac{\kappa}{\lambda+n+1}\right) \frac{\kappa^{2}j} {\lambda ( \lambda+n )
}\sum\limits_{s=0}^{n+1}\frac{(1-j) _{s} ( j-n ) _{s}(1-\kappa) _{s} (
1+\kappa ) _{s}}{s!(2) _{s} ( \lambda+1 ) _{s+1} ( -\lambda-n+1 ) _{s}}
\\
\qquad{}
=\left( \frac{\kappa}{\lambda+n+1}\right) a_{n,n+1-j},
\end{gather*}
where the index of summation is changed $s=i-1$.
\end{proof}

We consider the dif\/ferentiation formulae.
Set $\lambda_{1}:=\frac{m-\ell} {m}$, $\lambda_{2}:=\frac{\ell}{m}$.
In the previous section we used $\delta:=\frac{1}{2}-\frac{\ell}{m}$; thus $\delta=\frac{1}{2}-\lambda
_{2}=\lambda_{1}-\frac{1}{2}$.
For $n\geq0$ we have $\deg P_{n}^{(1) }=m(n+1) -\ell+1=m ( n+\lambda_{1} ) +1$ and
$\deg P_{n}^{(2) }=mn+\ell+1=m ( n+\lambda_{2} ) +1$, then (using $P_{-1}^{(1)
}=z^{1-\ell}t$ and $P_{-1} ^{(2) }=z^{1+\ell-m}\overline{t}$) for $j=1,2${\samepage
\begin{gather*}
\mathcal{D}z^{r}P_{n}^{(j) }= ( m(n+\lambda _{j}) +r+1 ) z^{r-1}P_{n}^{ (j) },\qquad 1\leq r\leq m-1,
\\
\mathcal{D}P_{n}^{(j) }= ( m(n+\lambda _{j}) +1 ) \left( z^{m-1}P_{n-1}^{(j) } +\frac{\kappa}{n+\lambda_{j}}\overline{z}^{m-1}\sigma_{0}P_{n-1}^{(j) }\right),
\\
\mathcal{D}^{2}P_{n}^{(j) }= ( m(n+\lambda _{j}) +1 ) m ( n+\lambda_{j} )
\left\{ 1-\left( \frac{\kappa}{n+\lambda_{j}}\right) ^{2}\right\} z^{m-2}P_{n-1}^{(j) },
\\
\overline{\mathcal{D}}z^{r}P_{n}^{(j) }=0,\qquad 0\leq r\leq m-2,
\\
\overline{\mathcal{D}}z^{m-1}P_{n}^{(j) }=-m\kappa \overline{z}^{m-2}\sigma_{0}P_{n}^{(j) }.
\end{gather*}
The other dif\/ferentiations follow from $\overline{\mathcal{D}}=\sigma _{0}\mathcal{D}\sigma_{0}$.}

Notice we have an evaluation formula for $p_{n}^{(j) }(1) $ ($j=1,2$): suppose $n=m (
k+\lambda_{j} ) +r$ with $1\leq r\leq m$ then
\begin{gather*}
p_{n}^{(1) }(1)=\frac{1}{2}\frac{ ( \lambda_{1}+\kappa ) _{k+1}}{ (
\lambda_{1} ) _{k+1}} ( t+\overline{t} ) +\frac{1}{2}\frac{ ( \lambda_{1}-\kappa ) _{k+1}}{ (
\lambda_{1} ) _{k+1}} ( t-\overline{t} ),
\\
p_{n}^{(2) }(1)=\frac{1}{2}\frac{ ( \lambda_{2}+\kappa ) _{k+1}}{ (
\lambda_{2} ) _{k+1}} ( t+\overline{t} ) -\frac{1}{2}\frac{ ( \lambda_{2}-\kappa ) _{k+1}}{ (
\lambda_{2} ) _{k+1}} ( t-\overline{t} ).
\end{gather*}

\section{The bilinear form on harmonic polynomials}  \label{section3}

Recall (formula~\eqref{pzform}) that the form is given by
\begin{gather*}
\big\langle p_{1} ( z,\overline{z} ) t+p_{2} ( z,\overline {z} ) \overline{t},q (
z,\overline{z},t,\overline{t} ) \big\rangle =\big\langle t,p_{1}^{\ast} ( 2\overline{\mathcal{D}
},2\mathcal{D} ) q\big\rangle |_{z=0}+\big\langle \overline{t},p_{2}^{\ast} (
2\overline{\mathcal{D}},2\mathcal{D} ) q\big\rangle |_{z=0}.
\end{gather*}
We will compute values for the harmonic polynomials.
Suppose $q$ is harmonic of degree $n$ and $p=\sum\limits_{j=0}^{n}a_{j}z^{n-j}\overline{z}^{j}t+\sum\limits
_{j=0}^{n}b_{j}z^{n-j}\overline{z}^{j}\overline{t}$ with arbitrary coef\/f\/icients $ \{ a_{j} \} $, $ \{
b_{j} \} $, then
\begin{gather}
2^{-n} \langle p,q \rangle=\left\langle t,\sum\limits_{j=0}
^{n}\overline{a_{j}}\mathcal{D}^{n-j}\overline{\mathcal{D}}^{j}q\right\rangle +\left\langle
\overline{t},\sum\limits_{j=0}^{n}\overline{b_{j}}\mathcal{D} ^{n-j}\overline{\mathcal{D}}^{j}q\right\rangle
\nonumber
\\
\hphantom{2^{-n} \langle p,q \rangle}{}
=\big\langle t,\big( \overline{a_{0}}\mathcal{D}^{n}+\overline{a_{n} }\overline{\mathcal{D}}^{n}\big)
q\big\rangle +\big\langle \overline {t},\big( \overline{b_{0}}\mathcal{D}^{n}+\overline{b_{n}}\overline
{\mathcal{D}}^{n}\big) q\big\rangle,\label{pDqz}
\end{gather}
because $\mathcal{D}^{n-j}\overline{\mathcal{D}}^{j}q=0$ for $1\leq j\leq n-1$.
Thus $\big\langle p_{n}^{(1) },q\big\rangle =2^{n}\big\langle t,\mathcal{D}^{n}q\big\rangle $ and
$\big\langle p_{n}^{(2) },q\big\rangle =2^{n}\big\langle \overline {t},\mathcal{D}^{n}q\big\rangle $;
this follows from the leading terms of $p_{n}^{(1) }$, $p_{n}^{(2) }$.

Def\/ine $\alpha_{n}$ by $\mathcal{D}^{n}p_{n}^{(1) }=\alpha_{n}t$ (we will show that there is no
$\overline{t}$ component).
If $1\leq n\leq m-\ell$ then $p_{n}^{(1) }=z^{n}t$ and $\mathcal{D}^{n} z^{n}t=n!t$ so $\alpha_{n}=n!$.
Suppose $n=m ( k+\lambda_{j} ) +r+1$ for some $k\geq0$ and $1\leq r\leq m-1$, then $p_{n}^{(1)
}=z^{r}P_{k}^{(1) }$ and $\mathcal{D}p_{n}^{(1) }= ( m ( n+\lambda_{j} )
+r+1 ) z^{r-1}P_{k}^{(1) }=nz^{r-1}P_{k}^{(1) }$.
Repeat this procedure to show $\mathcal{D}^{r}p_{n}^{(1) }=n ( n-1 ) \cdots ( n-r+1 )
P_{k}^{(1) }$, that is, $\alpha_{n}= ( -1 ) ^{r} ( -n ) _{r}\alpha_{n-r}$.
Now suppose $n=m ( k+\lambda_{j} ) +1$ so that $p_{n}^{(1) }=P_{k}^{(1) }$ and
\begin{gather}
\mathcal{D}^{2}P_{k}^{(1) }= ( m ( k+\lambda _{j} ) +1 ) m ( k+\lambda_{j} )
\left\{ 1-\left( \frac{\kappa}{k+\lambda_{1}}\right) ^{2}\right\} z^{m-2}P_{k-1}^{(1) }
\nonumber\\
\hphantom{\mathcal{D}^{2}P_{k}^{(1) }}{}
=n ( n-1 ) \left\{ 1-\left( \frac{\kappa}{k+\lambda_{1} }\right) ^{2}\right\} z^{m-2}P_{k-1}^{(1)
},\nonumber
\\
\mathcal{D}^{m}P_{k}^{(1) }=n ( n-1 ) \left\{ 1-\left( \frac{\kappa}{k+\lambda_{1}}\right)
^{2}\right\} \mathcal{D} ^{m-2}z^{m-2}P_{k-1}^{(1) }
\nonumber
\\
\hphantom{\mathcal{D}^{m}P_{k}^{(1) }}{}
=n(n-1) \left\{ 1-\left( \frac{\kappa}{k+\lambda_{1} }\right) ^{2}\right\} (n-2) \cdots (
n-m+1 ) P_{k-1}^{(1) }
\nonumber
\\
\hphantom{\mathcal{D}^{m}P_{k}^{(1) }}{}
=\frac{n!}{ ( n-m ) !}\frac{ ( \lambda_{1}-\kappa+k )  ( \lambda_{1}+\kappa+k ) }{ (
\lambda_{1}+k ) ^{2} }P_{k-1}^{(1) }.\label{DmPk}
\end{gather}
Induction on these results is used in the following

\begin{proposition}
For $k=0,1,2,\ldots$ and $n=m\left( k+\lambda_{1}\right) +1=m\left( k+1\right) -\ell+1$
\begin{gather*}
\mathcal{D}^{n}P_{k}^{(1) }=n!\frac{ ( \lambda _{1}-\kappa ) _{k+1} (
\lambda_{1}+\kappa ) _{k+1}}{ ( \lambda_{1} ) _{k+1}^{2}}t,
\\
\big\langle P_{k}^{(1) },P_{k}^{(1) }\big\rangle=2^{n+1}n!\frac{ (
\lambda_{1}-\kappa ) _{k+1} ( \lambda_{1}+\kappa ) _{k+1}}{ ( \lambda_{1} ) _{k+1}^{2}}.
\end{gather*}
\end{proposition}

\begin{proof}
Suppose $k=0$ then
\begin{gather*}
\mathcal{D}^{2}P_{0}^{(1) }= ( m-\ell+1 ) (m-\ell) \frac{ (
\lambda_{1}-\kappa )  ( \lambda _{1}+\kappa ) }{\lambda_{1}^{2}}z^{m-\ell-1}t
\end{gather*}
and $\mathcal{D}^{m-\ell-1}z^{m-\ell-1}t= ( m-\ell-1 ) !t$, so the formula is valid for $k=0$.
Suppose the formula is valid for some $k\geq0$ with $n=m ( k+1 ) -\ell+1$ then by~\eqref{DmPk}
\begin{gather*}
\mathcal{D}^{m}P_{k+1}^{(1) }=\frac{ ( n+m ) !}{n!}\frac{ ( \lambda_{1}-\kappa+k+1 )
 ( \lambda_{1} +\kappa+k+1 ) }{ ( \lambda_{1}+k+1 ) ^{2}}P_{k}^{(1) },
\\
\mathcal{D}^{n+m}P_{k+1}^{(1) }= ( n+m ) !\frac{\ ( \lambda_{1}-\kappa ) _{k+2} (
\lambda_{1} +\kappa ) _{k+2}}{ ( \lambda_{1} ) _{k+2}^{2}}t.
\end{gather*}
This completes the induction.
By Proposition~\ref{PzQ} the coef\/f\/icients of $z^{n}$ and $\overline{z}^{n}$ in $P_{k}^{(1) }$ are~$t$ and~$0$ respectively.
By equation~\eqref{pDqz} $\big\langle P_{k}^{(1) },P_{k}^{(1) }\big\rangle =2^{n}\big\langle
t,\mathcal{D} ^{n}P_{k}^{(1) }\big\rangle $.
The fact that $ \langle t,t \rangle =2$ f\/inishes the proof.
\end{proof}

\begin{proposition}
For $k=0,1,2,\ldots$ and $n=m ( k+\lambda_{2} ) +1=mk+\ell+1$
\begin{gather*}
\mathcal{D}^{n}P_{k}^{(2) }=n!\frac{ ( \lambda _{2}-\kappa ) _{k+1} (
\lambda_{2}+\kappa ) _{k+1}}{ ( \lambda_{2} ) _{k+1}^{2}}\overline{t},
\\
\big\langle P_{k}^{(2) },P_{k}^{(2) }\big\rangle=2^{n+1}n!\frac{ (
\lambda_{2}-\kappa ) _{k+1} ( \lambda_{2}+\kappa ) _{k+1}}{ ( \lambda_{2} ) _{k+1}^{2}}.
\end{gather*}
\end{proposition}

\begin{proof}
The same argument used in the previous propositions works here.
\end{proof}

The other basis polynomials can be handled as a~consequence.

\begin{proposition}
Suppose $j=1,2$, $k=0,1,2,\ldots$, $1\leq r\leq m-1$ and $n=m\left( k+\lambda_{j}\right) +r+1$ then $p_{n}^{\left(
j\right) }=z^{r} P_{k}^{(j) }$ and
\begin{gather*}
\big\langle p_{n}^{(j) },p_{n}^{(j) }\big\rangle =2^{n+1}n!\frac{ (
\lambda_{j}-\kappa ) _{k+1} ( \lambda_{j}+\kappa ) _{k+1}}{ ( \lambda_{j} ) _{k+1}^{2}}.
\end{gather*}
\end{proposition}

\begin{proof}
By Theorem~\ref{p2zp} $\mathcal{D}^{r}p_{n}^{(j) }=\frac {n!}{ ( n-r ) !}P_{k}^{(j) }$
and $\deg P_{k}^{(j) }=n-r=m ( k+\lambda_{j} ) +1$.
The preceding formulae for $\big\langle P_{k}^{(j) },P_{k}^{(j) }\big\rangle $ imply the
result.
\end{proof}

The trivial cases come from $\big\{ p_{n}^{(1) }:0\leq n\leq m-\ell\big\} $ and $\big\{ p_{n}^{ (
2 ) }:0\leq n\leq \ell\big\} $.
By Theorem~\ref{p2zp} $\mathcal{D}^{n}p_{n}^{(1) }=\mathcal{D}^{n}z^{n}t=n!t$ and $\big\langle
p_{n}^{(1) },p_{n}^{(1) }\big\rangle =2^{n+1}n!$ for $0\leq n\leq m-\ell$.
Similarly $\mathcal{D}^{n}p_{n}^{(2) }=\mathcal{D}^{n}z^{n}\overline{t}=n!\overline{t}$ and $\big\langle
p_{n}^{(2) },p_{n}^{(2) }\big\rangle =2^{n+1}n!$ for $0\leq n\leq\ell$.

\subsection{Orthogonal basis for the harmonics}

At each degree $\geq1$ we have the basis $\big\{ p_{n}^{(1) },\sigma_{0}p_{n}^{(1)
},p_{n}^{(2) },\sigma _{0}p_{n}^{(2) }:n\geq1\big\} $.
These polynomials are pairwise orthogonal with certain exceptions.
The leading terms and relevant properties are (the constants~$c_{n}$ are not specif\/ied, and may vary from line to line):{\samepage
\begin{enumerate}\itemsep=0pt
\item[1)]
$p_{n}^{(1) }$: $z^{n}t$, $\mathcal{D}^{n}p_{n}^{(1) }=c_{n}t$, $\overline{\mathcal{D}}p_{n}^{(1) }=0$ if $n+\ell\neq0\operatorname{mod}m$;

\item[2)]
$\sigma_{0}p_{n}^{(1) }$: $\overline{z}^{n}\overline {t}$, $\overline{\mathcal{D}}^{n}\sigma_{0}p_{n}^{(1) } =c_{n}\overline{t}$, $\mathcal{D}\sigma_{0}p_{n}^{(1) }=0$ if $n+\ell\neq0\operatorname{mod}m$;

\item[3)]
$p_{n}^{(2) }$: $z^{n}\overline{t}$, $\mathcal{D}^{n} p_{n}^{(2)
}=c_{n}\overline{t}$, $\overline{\mathcal{D}} p_{n}^{(2) }=0$ if $n-\ell\neq0\operatorname{mod}m$;

\item[4)]
$\sigma_{0}p_{n}^{(2) }$: $\overline{z}^{n}t$, $\overline {\mathcal{D}}^{n}\sigma_{0}p_{n}^{(2)
}=c_{n}t$, $\mathcal{D} \sigma_{0}p_{n}^{(2) }=0$ if $n-\ell\neq0\operatorname{mod}m$.
\end{enumerate}}

If $n\geq1$ then $\mathrm{span}\big\{ p_{n}^{(1) },\sigma _{0}p_{n}^{(1) }\big\}
\bot\, \mathrm{span}\big\{ p_{n}^{(2) },\sigma_{0}p_{n}^{(2) }\big\} $.
This is clear from the above properties when $n\pm\ell\neq0\operatorname{mod} m$ (note $\big\langle p_{n}^{(1) },p_{n}^{(2) }\big\rangle =0$ because $ \langle t,\overline{t} \rangle =0$).
One example suf\/f\/ices to demonstrate the other cases: suppose $n+\ell =0\operatorname{mod}m$ then
$n-\ell\neq0\operatorname{mod}m$ and $\big\langle p_{n}^{(1) },\sigma_{0}p_{n}^{(2)
}\big\rangle =2^{n}\big\langle t,\mathcal{D}^{n}\sigma_{0}p_{n}^{(2) }\big\rangle =0$.

Consider the special cases $n\pm\ell=0\operatorname{mod}m$.
Recall
\begin{gather*}
\overline{\mathcal{D}}z^{m-1}P_{k}^{(j) }=-m\kappa\sigma _{0}z^{m-2}P_{k}^{(j) }.
\end{gather*}
Using the leading term of $\sigma_{0}z^{m-1}P_{k}^{(j) }$ we f\/ind that
\begin{gather*}
\big\langle
\sigma_{0}z^{m-1}P_{k}^{(j) },z^{m-1}P_{k}^{(j) }\big\rangle = -2m\kappa\big\langle
z^{m-2}P_{k}^{(j) },z^{m-2}P_{k}^{(j) }\big\rangle \neq0.
\end{gather*}
In this case there are two dif\/ferent natural orthogonal bases.
One is $\big\{ p_{n}^{(j) }+\sigma_{0}p_{n}^{(j) },p_{n}^{(j)
}-\sigma_{0}p_{n}^{(j) }\big\} $.
First suppose $j=1$ and $n=m ( k+2 ) -\ell$ with $k\geq0$, then $p_{n}^{(1) }=z^{m-1}P_{k}^{(1) }$ and
\begin{gather*}
\big\langle p_{n}^{(1) },p_{n}^{(1) }\big\rangle=\big\langle \sigma_{0}p_{n}^{(1) },\sigma_{0}p_{n}^{(1) }\big\rangle =2^{n+1}n!\frac{( \lambda_{1}-\kappa) _{k+1}(
\lambda_{1}+\kappa) _{k+1} }{( \lambda_{1}) _{k+1}^{2}},
\\
\big\langle \sigma_{0}p_{n}^{(1) },p_{n}^{(1) }\big\rangle=2^{n}\big\langle
\overline{t},\overline{\mathcal{D}} ^{n}p_{n}^{(1) }\big\rangle =-2^{n}m\kappa\big\langle
\overline{t},\overline{\mathcal{D}}^{n-1}\sigma_{0}z^{m-2}P_{k}^{(1) }\big\rangle
\\
\hphantom{\big\langle \sigma_{0}p_{n}^{(1) },p_{n}^{(1) }\big\rangle}{}
=-2^{n}m\kappa(n-1) !\frac{ ( \lambda_{1}-\kappa ) _{k+1} ( \lambda_{1}+\kappa )
_{k+1}}{ ( \lambda_{1}\ ) _{k+1}^{2}} \langle \overline{t},\overline{t} \rangle.
\end{gather*}
These imply
\begin{gather*}
\big\langle p_{n}^{(1) }+\sigma_{0}p_{n}^{(1) },p_{n}^{(1)
}-\sigma_{0}p_{n}^{(1) } \big\rangle=\big\langle p_{n}^{(1) },p_{n}^{(1)
}\big\rangle -\big\langle \sigma_{0}p_{n}^{(1) },\sigma_{0}p_{n}^{(1) }\big\rangle =0,
\\
\big\langle p_{n}^{(1) }+\sigma_{0}p_{n}^{(1) },p_{n}^{(1)
}+\sigma_{0}p_{n}^{(1) } \big\rangle=2\big\langle p_{n}^{(1) },p_{n}^{(1)
}\big\rangle +2\big\langle \sigma_{0}p_{n}^{(1) },\sigma_{0}p_{n}^{(1) }\big\rangle
\\
\hphantom{\big\langle p_{n}^{(1) }+\sigma_{0}p_{n}^{(1) },p_{n}^{(1)
}+\sigma_{0}p_{n}^{(1) } \big\rangle}{}
=2^{n+2}n!\frac{ ( \lambda_{1}-\kappa ) _{k+2} ( \lambda_{1}+\kappa ) _{k+1}}{ (
\lambda_{1} ) _{k+2} ( \lambda_{1} ) _{k+1}},
\\
\big\langle p_{n}^{(1) }-\sigma_{0}p_{n}^{(1) },p_{n}^{(j)
}-\sigma_{0}p_{n}^{(1) } \big\rangle=2\big\langle p_{n}^{(1) },p_{n}^{(1)
}\big\rangle -2\big\langle \sigma_{0}p_{n}^{(1) },\sigma_{0}p_{n}^{(1) }\big\rangle
\\
\hphantom{\big\langle p_{n}^{(1) }-\sigma_{0}p_{n}^{(1) },p_{n}^{(j)
}-\sigma_{0}p_{n}^{(1) } \big\rangle}{}
=2^{n+2}n!\frac{ ( \lambda_{1}-\kappa ) _{k+1} ( \lambda_{1}+\kappa ) _{k+2}}{ (
\lambda_{1} ) _{k+1} ( \lambda_{1} ) _{k+2}}.
\end{gather*}
In these calculations we used
\begin{gather*}
n!-m\kappa(n-1) !=n!\left( \frac{n-m\kappa}{n}\right) =n!\left(
\frac{m ( k+1-\kappa ) +m-\ell}{m ( k+1 ) +m-\ell}\right) =n!\frac{\lambda_{1}
-\kappa+k+1}{\lambda_{1}+k+1}
\end{gather*}
 and
\begin{gather*}
n!+m\kappa(n-1) !=n!\frac{\lambda_{1}+\kappa+k+1}{\lambda_{1}+k+1}.
\end{gather*}
The result also applies when $k=-1$, that is, $n=m-\ell$ and $p_{n}^{(1) }=z^{m-\ell}t$.

Now suppose $j=2$ and $n=m( k+1) +\ell$ with $k\geq0$, then $p_{n}^{(2) }=z^{m-1}P_{k}^{(2) }$ and by a~similar argument
\begin{gather*}
\big\langle p_{n}^{(2) }+\sigma_{0}p_{n}^{(2) },p_{n}^{(2)
}+\sigma_{0}p_{n}^{(2) } \big\rangle=2^{n+2}n!\frac{ ( \lambda_{2}-\kappa ) _{k+2} (
\lambda_{2}+\kappa ) _{k+1}}{ ( \lambda_{2} ) _{k+2} ( \lambda_{2} ) _{k+1}},
\\
\big\langle p_{n}^{(2) }-\sigma_{0}p_{n}^{(2) },p_{n}^{(2)
}-\sigma_{0}p_{n}^{(2) } \big\rangle=2^{n+2}n!\frac{ ( \lambda_{2}-\kappa ) _{k+1} (
\lambda_{2}+\kappa ) _{k+2}}{ ( \lambda_{2} ) _{k+1} ( \lambda_{2} ) _{k+2}}.
\end{gather*}
Note $n=m ( k+1+\lambda_{2} ) $ and $n\pm m\kappa=m ( \lambda_{2}\pm\kappa+k+1 ) $.
The formula also applies to $n=\ell$ with $k=-1$ and $p_{n}^{(2) }=z^{\ell}\overline{t}$.

The other natural basis is $\big\{ p_{n}^{(j) }+\frac{m\kappa }{n}\sigma_{0}p_{n}^{(j)
},\sigma_{0}p_{n}^{(j) }\big\} $ (with $p_{n}^{(j) }=z^{m-1}P_{k}^{(j) }$); the
idea is to form a~linear combination $f_{n}^{(j) }=p_{n}^{(j) }+c\sigma_{0}p_{n}^{(j) }$ so that $\overline{\mathcal{D}}f=0$, which implies $\big\langle \sigma_{0} p_{n}^{(j)
},f_{n}^{(j) }\big\rangle =0$.
Indeed
\begin{gather*}
\overline{\mathcal{D}}f_{n}^{(j) }=-m\kappa\sigma_{0} z^{m-2}P_{k}^{(j)
}+cn\sigma_{0}z^{m-2}P_{k}^{(j) }.
\end{gather*}
Another expression for $f_{n}^{(j) }$ is $\frac{1}{z} P_{k+1}^{(j) }$, that is,
\begin{gather*}
f_{n}^{(1) }=z^{m-\ell}Q_{k}^{(1) }\big( \kappa,\lambda_{1};z^{m},\overline{z}^{m}\big)
t+\overline{z}^{m-\ell} Q_{k}^{(2) }\big( \kappa,\lambda_{1};z^{m},\overline{z} ^{m}\big)
\overline{t}, \qquad n=m ( k+1 ) -\ell,
\\
f_{n}^{(2) }=z^{\ell}Q_{n}^{(1) }\big( \kappa,\lambda_{2};z^{m},\overline{z}^{m}\big)
\overline{t}+\overline {z}^{\ell}Q_{n}^{(2) }\big( \kappa,\lambda_{2};z^{m},\overline{z}^{m}\big)
t,\qquad n=mk+\ell.
\end{gather*}
To compute $\big\langle f_{n}^{(j) },f_{n}^{(j) }\big\rangle $ (by use of
$\mathcal{D}\sigma_{0}=\sigma_{0}\overline {\mathcal{D}}$)
\begin{gather*}
\mathcal{D}^{m-1}f_{n}^{(j) }=\mathcal{D}^{m-2}\left( nz^{m-1}P_{k}^{(j)
}-m\kappa\frac{m\kappa}{n}z^{m-2} P_{k}^{(j) }\right)
\\
\hphantom{\mathcal{D}^{m-1}f_{n}^{(j) }}{}
=n\left( 1-\left( \frac{m\kappa}{n}\right) ^{2}\right) \mathcal{D} ^{m-2}z^{m-1}P_{k}^{(j) }
=\frac{n!}{ ( n-m+1 ) !}\left( 1-\left( \frac{m\kappa} {n}\right) ^{2}\right) P_{k}^{(j) }.
\end{gather*}
If $j=1$ then $n=m ( k+2 ) -\ell$ and
\begin{gather*}
1-\left( \frac{m\kappa} {n}\right) ^{2}=\frac{(
k+1+\lambda_{1}-\kappa) ( k+1+\lambda_{1}+\kappa) }{( k+1+\lambda_{1}) ^{2}},
\end{gather*} and by the
previous results (also valid for $k=-1$)
\begin{gather*}
\big\langle f_{n}^{(1) },f_{n}^{(1) }\big\rangle =2^{n+1}n!\frac{(
\lambda_{1}-\kappa) _{k+2}( \lambda_{1}+\kappa) _{k+2}}{( \lambda_{1}) _{k+2}^{2}}.
\end{gather*}
Similarly if $j=2$ then $n=m( k+1) +\ell$ and
\begin{gather*}
\big\langle f_{n}^{(2) },f_{n}^{(2) }\big\rangle =2^{n+1}n!\frac{(
\lambda_{2}-\kappa) _{k+2}( \lambda_{2}+\kappa) _{k+2}}{( \lambda_{2}) _{k+2}^{2}}.
\end{gather*}
Observe that $-\frac{\ell}{m}<\kappa<\frac{\ell}{m}=\lambda_{2}$ implies all the pairings $\langle p,p\rangle
>0$ (by def\/inition $\lambda _{2}<\lambda_{1}$).

\section{The Gaussian kernel}  \label{section4}

In the previous sections we showed that the forms $ \langle \cdot,\cdot \rangle $ and $ \langle
\cdot,\cdot \rangle _{G}$ are positive-def\/inite on the harmonic polynomials exactly when $-\frac{\ell}
{m}<\kappa<\frac{\ell}{m}$.
By use of general formulae (\cite[(4), p.~4]{Dunkl2013}) with $\gamma ( \kappa;\tau ) =0$) for
$\Delta_{\kappa}^{k} (  ( z\overline{z} ) ^{n}f ( z,\overline{z},t,\overline{t} )  ) $, where
$\Delta_{\kappa}f=0$ it follows that the forms are positive-def\/inite on~$\mathcal{P}_{m,\ell}$ for the same values of
$\kappa$.
This together with the property $ \langle pf,g \rangle _{G}= \langle f,p^{\ast}g \rangle _{G}$ for any
scalar polynomial $p ( z,\overline{z} ) $ suggests there should be an integral formulation of $ \langle
\cdot,\cdot \rangle _{G}$.
The construction begins with suitable generalizations of the integral for the scalar case and the group~$I_{2} (
m ) $, namely $\vert z^{m}-\overline{z}^{m} \vert ^{2\kappa}\exp\big( {-}\frac{z\overline{z} }{2}\big)
$.

In~\cite[Section~5]{Dunkl2013} we derived conditions, holding for any f\/inite ref\/lection group, on the matrix function~$K$ associated with
the Gaussian inner product.
The formula for $W$ is stated in~\eqref{defK}, and the general conditions specialize to:
\begin{enumerate}\itemsep=0pt
\item[1)]
$K ( zw ) =\tau_{\ell} ( w ) ^{-1}K ( z ) \tau_{\ell} ( w ),w\in W$;

\item[2)]
$K$ satisf\/ies a~boundary condition at the walls of the fundamental chamber  $\mathcal{C}:=\big\{
z=re^{\mathrm{i}\theta}:0<\theta<\frac{\pi} {m},r>0\big\} $; see equation~\eqref{Kbdry} below;

\item[3)]
$K=L^{\ast}ML$ where $M$ is a~constant positive-def\/inite matrix and $L$ satisf\/ies a~dif\/ferential system.
\end{enumerate}

The dif\/ferential system (for a~$2\times2$ matrix function~$L$) is
\begin{gather*}
\partial_{z}L ( z,\overline{z} )=\kappa L ( z,\overline{z} )
\sum\limits_{j=0}^{m-1}\frac{1}{z-\overline{z}\omega^{j}} \tau_{\ell} ( \sigma_{j} ),
\\
\partial_{\overline{z}}L ( z,\overline{z} )=\kappa L ( z,\overline{z} )
\sum\limits_{j=0}^{m-1}\frac{-\omega^{j}}{z-\overline{z} \omega^{j}}\tau_{\ell} ( \sigma_{j} ).
\end{gather*}
(As motivation for this formulation note that the analogous scalar weight function $ \vert
z^{m}-\overline{z}^{m} \vert ^{2\kappa}=\overline {h}h$ where $h$ satisf\/ies $\partial_{z}h=\kappa
h\sum\limits_{j=0}^{m-1}\frac {1}{z-\overline{z}\omega^{j}}$ and $\partial_{\overline{z}}h=\kappa
h\sum\limits_{j=0}^{m-1}\frac{-\omega^{j}}{z-\overline{z}\omega^{j}}$.) It follows that
\begin{gather*}
 ( z\partial_{z}+\overline{z}\partial_{\overline{z}} ) L=0,
\end{gather*}
thus $L$ is positively homogeneous of degree $0$ and depends only on $e^{\mathrm{i}\theta}$ for
$z=re^{\mathrm{i}\theta}$ ($r>0$, $-\pi<\theta\leq\pi $).
Further $\partial_{\theta}=\mathrm{i} ( z\partial_{z}-\overline {z}\partial_{\overline{z}} ) $ and thus
\begin{gather*}
\partial_{\theta}L=\mathrm{i}\kappa L\sum\limits_{j=0}^{m-1}\frac{z+\overline
{z}\omega^{j}}{z-\overline{z}\omega^{j}}\tau_{\ell} ( \sigma_{j} ).
\end{gather*}
We use the following to compute the sums involving $z-\overline{z}\omega^{j}$.

\begin{lemma}
\label{sum1}
For $n\in\mathbb{Z}$ let $r=n-1\operatorname{mod}m$ and $0\leq r<m$ then $($indeterminate $w)$
\begin{gather*}
\sum\limits_{j=0}^{m-1}\frac{\omega^{jn}}{w-\omega^{j}}=\frac{mw^{r}}{w^{m} -1},
\qquad
\sum\limits_{j=0}^{m-1}\frac{\omega^{jn}}{z-\overline{z}\omega^{j}}=\frac
{mz^{r}\overline{z}^{m-1-r}}{z^{m}-\overline{z}^{m}}.
\end{gather*}
\end{lemma}

\begin{proof}
Assume $ \vert w \vert <1$ then
\begin{gather*}
\sum\limits_{j=0}^{m-1}\frac{\omega^{jn}}{w-\omega^{j}}=-\sum\limits_{j=0}^{m-1}\frac {\omega^{j(n-1)
}}{1-w\omega^{-j}}=-\sum\limits_{s=0}^{\infty}w^{s} \sum\limits_{j=0}^{m-1}\omega^{j ( n-1-s ) }.
\end{gather*}
The inner sum equals $m$ if $n-1=s\operatorname{mod}m$, that is $s=r+km$ for $k\geq0$, otherwise the sum vanishes; thus
\begin{gather*}
\sum\limits_{j=0}^{m-1}\frac{\omega^{jn}}{w-\omega^{j}}=-m\sum\limits_{k=0}^{\infty} w^{r+km}=\frac{-mw^{r}}{1-w^{m}}.
\end{gather*}
The identity is valid for all $w\notin \{ \omega^{j} \} $.
Next
\begin{gather*}
\sum\limits_{j=0}^{m-1}\frac{\omega^{jn}}{z-\overline{z}\omega^{j}}=\frac
{1}{\overline{z}}\sum\limits_{j=0}^{m-1}\frac{\omega^{jn}}{\frac{z}{\overline{z}
}-\omega^{j}}=\frac{m}{\overline{z}}\frac{z^{r}\overline{z}^{-r}} {z^{m}\overline{z}^{-m}-1}
=\frac{mz^{r}\overline{z}^{m-1-r}}{z^{m}-\overline{z}^{m}}.\tag*{\qed}
\end{gather*}
\renewcommand{\qed}{}
\end{proof}

\begin{corollary}
For $1\leq\ell\leq m-1$
\begin{gather*}
\sum\limits_{j=0}^{m-1}\frac{z+\overline{z}\omega^{j}}{z-\overline{z}\omega^{j} }\omega^{j\ell}=2m\frac{z^{\ell}\overline{z}^{m-\ell}}{z^{m}-\overline {z}^{m}},
\qquad
\sum\limits_{j=0}^{m-1}\frac{z+\overline{z}\omega^{j}}{z-\overline{z}\omega^{j} }\omega^{-j\ell}=2m\frac{z^{m-\ell}\overline{z}^{\ell}}{z^{m} -\overline{z}^{m}}.
\end{gather*}
\end{corollary}

Write $L=L_{1}t+L_{2}\overline{t}$, then
\begin{gather*}
\partial_{\theta}L_{1} ( \theta )=\frac{m\kappa}{\sin m\theta }e^{\mathrm{i}\theta ( m-2\ell )
}L_{2} ( \theta ),
\qquad
\partial_{\theta}L_{2} ( \theta )=\frac{m\kappa}{\sin m\theta }e^{-\mathrm{i}\theta ( m-2\ell )
}L_{1} ( \theta ).
\end{gather*}
As in~\cite[Theorem~4.3]{Dunkl1993} let $\delta=\frac{1}{2}-\frac{\ell}{m}$, $\widetilde{L_{1}}=e^{-\mathrm{i}m\delta\theta}L_{1}$,
$\widetilde{L_{2} }=e^{\mathrm{i}m\delta\theta}L_{2}$, then
\begin{gather*}
\partial_{\theta}\widetilde{L}_{1} ( \theta )=-\mathrm{i} m\delta\widetilde{L}_{1} ( \theta )
+\frac{m\kappa}{\sin m\theta }\widetilde{L}_{2} ( \theta ),
\qquad
\partial_{\theta}\widetilde{L}_{2}(\theta)=\frac{m\kappa }{\sin m\theta}\widetilde{L}_{1} (\theta) +\mathrm{i} m\delta\widetilde{L}_{2}(\theta).
\end{gather*}
Now changing $ \{ t,\overline{t} \} $ to real coordinates $t=s_{1}+\mathrm{i}s_{2}$ we write
$\widetilde{L}_{1}t+\widetilde{L} _{2}\overline{t}=\big( \widetilde{L}_{1}+\widetilde{L}_{2}\big)
s_{1}+\mathrm{i}\big( \widetilde{L}_{1}-\widetilde{L}_{2}\big) s_{2}$.
Let $g_{1}=\widetilde{L}_{1}+\widetilde{L}_{2},g_{2}=\mathrm{i}\big( \widetilde{L}_{1}-\widetilde{L}_{2}\big) $ and
$\phi=m\theta$.
The system is transformed to
\begin{gather*}
\partial_{\phi}g_{1}=\frac{\kappa}{\sin\phi}g_{1}-\delta g_{2},
\qquad
\partial_{\phi}g_{2}=\delta g_{1}-\frac{\kappa}{\sin\phi}g_{2}.
\end{gather*}
Introduce the variable $v=\sin^{2}\frac{\phi}{2}$ and let
\begin{gather*}
g_{1}=\left( \frac{v}{1-v}\right) ^{\kappa/2}h_{1}(v), \qquad g_{2}=\left( \frac{v}{1-v}\right) ^{(\kappa+1) /2} h_{2}(v).
\end{gather*}
Here we restrict $\theta$ to the fundamental chamber $\big\{ \theta:0<\theta<\frac{\pi}{m}\big\} $ so that $0<v<1$
and $\sin\phi=2 ( v(1-v)  ) ^{1/2}$ and $\frac{\partial}{\partial\phi }v= ( v(1-v)
) ^{1/2}$. Then
\begin{gather*}
\partial_{v}h_{1}=-\frac{\delta}{1-v}h_{2},
\qquad
\partial_{v}h_{2}=\frac{\delta}{v}h_{1}-\frac{\kappa+\frac{1}{2}}{v(1-v) }h_{2},
\end{gather*}
from which follows
\begin{gather*}
v(1-v) \partial_{v}^{2}h_{1}+\left( \kappa+\frac{1} {2}-v\right) \partial_{v}h_{1}+\delta^{2}h_{1}=0.
\end{gather*}
This is the hypergeometric equation with solutions
\begin{gather*}
h_{1}^{(1) }={}_{2}F_{1}\left(   \genfrac{}{}{0pt}{}{\delta,-\delta}{\frac{1}{2}+\kappa} ;v\right),
\\
h_{1}^{(2) }=v^{\frac{1}{2}-\kappa}\;{} _{2}F_{1}\left(
\genfrac{}{}{0pt}{}{\frac{1}{2}-\kappa-\delta,\frac{1}{2}-\kappa+\delta }{\frac{3}{2}-\kappa} ;v\right)
=v^{\frac{1}{2}-\kappa}(1-v) ^{\frac{1}{2}+\kappa}\;{} _{2} F_{1}\left(
\genfrac{}{}{0pt}{}{1+\delta,1-\delta}{\frac{3}{2}-\kappa} ;v\right).
\end{gather*}
By use of the dif\/ferentiations
\begin{gather*}
\partial_{v}\;{} _{2}F_{1}\left(   \genfrac{}{}{0pt}{}{a,b}{c} ;v\right)=\frac{ab}{c}\;{}_{2}F_{1}\left(
\genfrac{}{}{0pt}{}{a+1,b+1}{c+1} ;v\right),
\\
\partial_{v}v^{c}\;{} _{2}F_{1}\left(   \genfrac{}{}{0pt}{}{a,b}{c+1} ;v\right)=cv^{c-1}\;{} _{2}F_{1}\left(
\genfrac{}{}{0pt}{}{a,b}{c} ;v\right),
\end{gather*}
we f\/ind
\begin{gather*}
h_{2}^{(1) }=\frac{\delta}{\frac{1}{2}+\kappa}(1-v)\;{} _{2}F_{1}\left(
\genfrac{}{}{0pt}{}{1+\delta,1-\delta}{\frac{3}{2}+\kappa} ;v\right),
\\
h_{2}^{(2) }=-\frac{\frac{1}{2}-\kappa}{\delta}(1-v) v^{-\frac{1}{2}-\kappa}\;{} _{2}F_{1}\left(
\genfrac{}{}{0pt}{}{\frac{1}{2}-\kappa-\delta,\frac{1}{2}-\kappa+\delta }{\frac{1}{2}-\kappa} ;v\right)
\\
\hphantom{h_{2}^{(2) }}{}
=-\frac{\frac{1}{2}-\kappa}{\delta}(1-v) ^{\kappa+\frac{1} {2}}v^{-\frac{1}{2}-\kappa}\;{}_{2}F_{1}\left(
\genfrac{}{}{0pt}{}{\delta,-\delta}{\frac{1}{2}-\kappa} ;v\right).
\end{gather*}
All these solutions can be expressed in terms of the following
\begin{gather}
f_{1} ( \kappa,\delta;v ) :=v^{\kappa/2}(1-v) ^{-\kappa/2}\;{}_{2}F_{1}\left(
\genfrac{}{}{0pt}{}{\delta,-\delta}{\frac{1}{2}+\kappa} ;v\right),
\nonumber\\
f_{2} ( \kappa,\delta;v ) :=\frac{\delta}{\frac{1}{2}+\kappa }v^{ ( \kappa+1 ) /2}(1-v)
^{(1-\kappa) /2}\;{}_{2}F_{1}\left(   \genfrac{}{}{0pt}{}{1+\delta,1-\delta}{\frac{3}{2}+\kappa} ;v\right).
\label{f1f2def}
\end{gather}
Then (by replacing $h_{j}^{(2) }$ by $\frac{-\delta}{\frac{1} {2}-\kappa}h_{j}^{(2) }$ for
$j=1,2$)
\begin{alignat*}{3}
& g_{1}^{(1) }=\left( \frac{v}{1-v}\right) ^{\kappa/2} h_{1}^{(1) }=f_{1} (
\kappa,\delta;v ),\qquad &&
g_{2}^{(1) }=\left( \frac{v}{1-v}\right) ^{\left( \kappa+1\right) /2}h_{2}^{(1) }=f_{2} (
\kappa,\delta;v ),&
\\
& g_{1}^{(2) }=\left( \frac{v}{1-v}\right) ^{\kappa/2} h_{1}^{(2) }=-f_{2} (
-\kappa,\delta;v ),\qquad &&
g_{2}^{(2) }=\left( \frac{v}{1-v}\right) ^{\left( \kappa+1\right) /2}h_{2}^{(2) }=f_{1} (
-\kappa,\delta;v ).&
\end{alignat*}
The Wronskian of the system is $f_{1} ( \kappa,\delta;v ) f_{1} ( -\kappa,\delta;v ) +f_{2} (
\kappa,\delta;v ) f_{2} ( -\kappa,\delta;v ) =1$.
Returning to~$L_{1}$,~$L_{2}$ we f\/ind (for $j=1,2$)
\begin{gather*}
L_{1}^{(j) }=e^{\mathrm{i}m\delta\theta}\widetilde{L} _{1}^{(j)
}=\frac{1}{2}e^{\mathrm{i}m\delta\theta}\big( g_{1}^{(j) }-\mathrm{i}g_{2}^{(j) }\big),
\qquad
L_{2}^{(j) }=e^{-\mathrm{i}m\delta\theta}\widetilde{L} _{2}^{(j)
}=\frac{1}{2}e^{-\mathrm{i}m\delta\theta}\big( g_{1}^{(j) }+\mathrm{i}g_{2}^{(j) }\big),
\end{gather*}
and in the real coordinate system $t=s_{1}+\mathrm{i}s_{2}$
\begin{gather*}
L_{j1} :=L_{1}^{(j) }+L_{2}^{(j) }=\cos m\delta\theta~g_{1}^{(j) }+\sin
m\delta\theta~g_{2}^{(j) },
\\
L_{j2} :=\mathrm{i}\big( L_{1}^{(j) }-L_{2}^{(j) }\big) =-\sin m\delta\theta~g_{1}^{(j) }+\cos m\delta\theta~g_{2}^{(j) }.
\end{gather*}
So a~fundamental solution in the real coordinates (for $0<\phi=m\theta<\pi$, recall $v=\sin^{2}\frac{\phi}{2}$) is
\begin{gather}
L ( \phi )=M_{f}(v) M_{\delta} ( \phi ),
\label{Mdef}
\\
M_{f}(v)=\left[
\begin{matrix}
f_{1} ( \kappa,\delta;v ) & f_{2} ( \kappa,\delta;v )
\\
-f_{2} ( -\kappa,\delta;v ) & f_{1} ( -\kappa,\delta;v )
\end{matrix}
\right], \qquad M_{\delta} ( \phi ) =\left[
\begin{matrix}
\cos\delta\phi & -\sin\delta\phi
\\
\sin\delta\phi & \cos\delta\phi
\end{matrix}
\right].
\nonumber
\end{gather}
(In the trivial case $\kappa=0$ it follows from standard hypergeometric series formulae that $f_{1} (
0,\delta;v )$ $=\cos\delta\phi$ and $f_{2} ( 0,\delta;v ) =\sin\delta\phi$, see Lemma~\ref{f12k0}.) As
described in~\cite{Dunkl2013} the kernel for the Gaussian inner product in the fundamental chamber $\mathcal{C}$ is
given by $K(\phi) e^{- ( x_{1}^{2}+x_{2}^{2} ) /2}$ where
\begin{gather*}
K(\phi) =L(\phi) ^{T}ML(\phi)
\end{gather*}
and $M$ is a~constant positive-def\/inite matrix determined by boundary conditions at the walls, $\theta=0$ and
$\theta=\frac{\pi}{m}$ in this case.
Suppose $\sigma$ corresponds to one of the walls (that is, the hyperplane f\/ixed by $\sigma$) and $x_{0}\neq0$ is
a~boundary point of $\mathcal{C}$ with $x_{0}\sigma=x_{0}$ then the condition is: for any vectors $\xi,\eta$ in the
module for $\tau_{\ell}$ (the underlying representation) if $\xi\tau_{\ell }(\sigma) =\xi$ and
$\eta\tau_{\ell}(\sigma) =-\eta$ then
\begin{gather*}
\lim_{x\rightarrow x_{0},\, x\in\mathcal{C}}\xi K(x) \eta^{T}=0,
\end{gather*}
equivalently
\begin{gather}
\lim_{x\rightarrow x_{0},\, x\in\mathcal{C}}\big( K(x) -\tau_{\ell}(\sigma) K(x)
\tau_{\ell}(\sigma) \big) =0.
\label{Kbdry}
\end{gather}
By the invariance property of $K$ the condition is equivalent to
\begin{gather*}
\lim\limits_{x\rightarrow0,\, x_{0}+x\in\mathcal{C}}\big( K ( x_{0}+x ) -K ( x_{0}+x\sigma ) \big) =0,
\end{gather*}
a~weak type of continuity.
Near $\phi=0$ consider $\tau_{\ell} ( \sigma_{0} ) =\left[
\begin{matrix}
1 & 0
\\
0 & -1
\end{matrix}
\right] $ and the boundary condition becomes
\begin{gather*}
\lim_{\phi\rightarrow0_{+}}K(\phi) _{12}=0.
\end{gather*}
Let $M_{0}=\left[
\begin{matrix}
a_{1} & a_{2}
\\
a_{2} & a_{3}
\end{matrix}
\right] $ and replace $M_{f}(v) $ by (with nonzero constants $c_{1}$, $c_{2}$ depending on~$\delta$ and~$\kappa$, see~\eqref{f1f2def})
\begin{gather*}
L^{\prime}(\phi):=\left[
\begin{matrix}
2^{-\kappa}\phi^{\kappa} & c_{1}\phi^{\kappa+1}
\\
c_{2}\phi^{1-\kappa} & 2^{\kappa}\phi^{-\kappa}
\end{matrix}
\right],
\end{gather*}
then $\lim\limits_{\phi\rightarrow0}M_{\delta}(\phi) =I$ and
\begin{gather*}
\big( L^{\prime}(\phi) ^{T}M_{0}L^{\prime}(\phi) \big)
_{12}=2^{-\kappa}c_{2}a_{1}\phi^{1+2\kappa}-2^{\kappa}c_{1}a_{3} \phi^{1-2\kappa}+a_{2}\big(
1-c_{1}c_{2}\phi^{2}\big).
\end{gather*}
The boundary condition immediately implies $a_{2}=0$, then the positive-def\/inite condition implies $a_{1}a_{3}\neq0$ and
thus $1\pm2\kappa >0$.
This bound ($-\frac{1}{2}<\kappa<\frac{1}{2}$) also implies the integrability of $K$ on the circle (and with respect to
$e^{- \vert x \vert ^{2}/2}dm_{2}(x) $).

Next we consider the behavior of $L(\phi) $ near $\phi=\pi$, that is, $v$ near $1$.
We use (see~\cite[(15.8.4)]{Olver2010})
\begin{gather*}
{}_{2}F_{1}\left(   \genfrac{}{}{0pt}{}{a,b}{c} ;v\right) =\frac{\Gamma(c) \Gamma(c-a-b) }
{\Gamma( c-a) \Gamma( c-b) }\;{}_{2}F_{1}\left(   \genfrac{}{}{0pt}{}{a,b}{1+a+b-c} ;1-v\right)
\\
\hphantom{{}_{2}F_{1}\left(   \genfrac{}{}{0pt}{}{a,b}{c} ;v\right) =}{}
+(1-v) ^{c-a-b}\frac{\Gamma(c) \Gamma ( a+b-c ) }{\Gamma ( a ) \Gamma (
b ) }\;{}_{2} F_{1}\left(   \genfrac{}{}{0pt}{}{c-a,c-b}{1-a-b+c} ;1-v\right).
\end{gather*}
Def\/ine
\begin{gather*}
H ( \kappa,\delta ):=\frac{\Gamma\left( \frac{1}{2} +\kappa\right) ^{2}}{\Gamma\left(
\frac{1}{2}+\kappa+\delta\right) \Gamma\left( \frac{1}{2}+\kappa-\delta\right) }.
\end{gather*}
The formula $\Gamma ( u ) \Gamma ( 1-u ) =\frac{\pi }{\sin\pi u}$ implies
\begin{gather*}
\frac{\Gamma\left( \frac{1}{2}+\kappa\right) \Gamma\left( -\frac{1} {2}-\kappa\right) }{\Gamma(\delta)
\Gamma(-\delta) }=\left( \frac{\delta}{\frac{1}{2}+\kappa}\right) \frac{\sin \pi\delta}{\cos\pi\kappa}.
\end{gather*}
Then
\begin{gather}
f_{1}(\kappa,\delta;v)=H(\kappa,\delta) f_{1}(-\kappa,\delta;1-v)
+\frac{\sin\pi\delta}{\cos\pi\kappa }f_{2}(\kappa,\delta;1-v),
\nonumber
\\
f_{2}(\kappa,\delta;v)=\frac{\sin\pi\delta}{\cos\pi\kappa }f_{1}(\kappa,\delta;1-v) -H (\kappa,\delta) f_{2}(-\kappa,\delta;1-v).
\label{fv21mv}
\end{gather}
Note $H(\kappa,\delta) H ( -\kappa,\delta ) +\left( \frac{\sin\pi\delta}{\cos\pi\kappa}\right)
^{2}=1$.
Summarize the transformation rules as
\begin{gather*}
\left[
\begin{matrix}
f_{1}(\kappa,\delta;v) & f_{2}(\kappa,\delta;v)
\\
-f_{2}(-\kappa,\delta;v)\!\! & f_{1}(-\kappa,\delta;v)
\end{matrix}
\right] =
\left[
\begin{matrix}
H(\kappa,\delta)\! & \frac{\sin\pi\delta}{\cos\pi\kappa}
\\
-\frac{\sin\pi\delta}{\cos\pi\kappa}\! & H ( -\kappa,\delta )
\end{matrix}
\right]\! \left[
\begin{matrix}
f_{1}(-\kappa,\delta;1-v)\!\! & -f_{2}(-\kappa,\delta;1-v)
\\
f_{2}(\kappa,\delta;1-v) & f_{1}(\kappa,\delta;1-v)
\end{matrix}
\right],
\end{gather*}
or set
\begin{gather*}
M_{H}=\left[
\begin{matrix}
H(\kappa,\delta) & \frac{\sin\pi\delta}{\cos\pi\kappa}
\\
-\frac{\sin\pi\delta}{\cos\pi\kappa} & H(-\kappa,\delta)
\end{matrix}
\right],\qquad M_{\sigma}=\left[
\begin{matrix}
0 & 1
\\
1 & 0
\end{matrix}
\right],
\end{gather*}
then
\begin{gather*}
M_{f}(v) =M_{H}M_{\sigma}M_{f}(1-v) M_{\sigma},
\end{gather*}
the same relation holds with $v$ and $1-v$ interchanged (indeed $ ( M_{H}M_{\sigma} ) ^{2}=I$).
The ref\/lection corresponding to the wall $\theta=\frac{\pi}{m}$ is $\sigma_{1}$ ($z\sigma_{1}=\overline{z}\omega$ so
$e^{\mathrm{i}\pi/m}\sigma_{1}=e^{\mathrm{i}\pi/m}$).
Thus the boundary condition is
\begin{gather*}
\lim_{\phi\rightarrow\pi_{-}}\big( K(\phi) -\tau_{\ell} ( \sigma_{1} ) K(\phi)
\tau_{\ell} ( \sigma_{1} ) \big) =0,
\end{gather*}
and (in the real coordinates)
\begin{gather*}
\tau_{\ell}\left( \sigma_{1}\right) =\left[
\begin{matrix}
\cos\frac{2\pi\ell}{m} & \sin\frac{2\pi\ell}{m}
\vspace{1mm}\\
\sin\frac{2\pi\ell}{m} & -\cos\frac{2\pi\ell}{m}
\end{matrix}
\right].
\end{gather*}
Then
\begin{gather*}
K(\phi) -\tau_{\ell}\left( \sigma_{1}\right) K(\phi) \tau_{\ell}\left( \sigma_{1}\right)
=q(\phi) \left[
\begin{matrix}
\sin\frac{2\pi\ell}{m} & -\cos\frac{2\pi\ell}{m}
\vspace{1mm}\\
-\cos\frac{2\pi\ell}{m} & -\sin\frac{2\pi\ell}{m}
\end{matrix}
\right],
\\
q(\phi) =\left\{ K(\phi) _{11}-K(\phi) _{22}\right\} \sin\frac{2\pi\ell}{m}-2K\left(
\phi\right) _{12}\cos\frac{2\pi\ell}{m}.
\end{gather*}
We want to determine $M_{0}$ so that $\lim\limits_{\phi\rightarrow\pi_{-} }q(\phi) =0$.
We have
\begin{gather*}
K(\phi) =M_{\delta}(\phi) ^{T}M_{f}(v) ^{T}M_{0}M_{f}(v)
M_{\delta}(\phi).
\end{gather*}
Let $\psi=\pi-\phi$ then the f\/irst-order approximations (as $\psi \rightarrow0_{+}$) are
\begin{gather*}
1-v \sim\frac{1}{4}\psi^{2},
\\
f_{1}(\kappa,\delta;1-v) \sim\left( \frac{\psi}{2}\right) ^{-\kappa},\qquad f_{2}(\kappa,\delta;1-v)
\sim\frac{\delta}{\frac {1}{2}+\kappa}\left( \frac{\psi}{2}\right) ^{1-\kappa},
\\
~f_{1}(-\kappa,\delta;1-v) \sim\left( \frac{\psi} {2}\right) ^{\kappa}, \qquad f_{2} (
-\kappa,\delta;1-v ) \sim \frac{\delta}{\frac{1}{2}-\kappa}\left( \frac{\psi}{2}\right) ^{1+\kappa},
\end{gather*}
and
\begin{gather*}
M_{f}(v) \sim M_{H}\left[
\begin{matrix}
\left( \frac{\psi}{2}\right) ^{\kappa} & -\frac{\delta}{\frac{1}{2}-\kappa }\left( \frac{\psi}{2}\right)^{1+\kappa}
\vspace{1mm}\\
\frac{\delta}{\frac{1}{2}+\kappa}\left( \frac{\psi}{2}\right) ^{1-\kappa} & \left( \frac{\psi}{2}\right)^{-\kappa}
\end{matrix}
\right].
\end{gather*}
We also f\/ind
\begin{gather*}
M_{\delta}\left( \pi-\psi\right) \sim M_{1}+\psi\delta M_{2},
\\
M_{1}=\left[
\begin{matrix}
\sin\frac{\pi\ell}{m} & -\cos\frac{\pi\ell}{m}
\vspace{1mm}\\
\cos\frac{\pi\ell}{m} & \sin\frac{\pi\ell}{m}
\end{matrix}
\right],\qquad M_{2}=\left[
\begin{matrix}
\cos\frac{\pi\ell}{m} & \sin\frac{\pi\ell}{m}
\vspace{1mm}\\
-\sin\frac{\pi\ell}{m} & \cos\frac{\pi\ell}{m}
\end{matrix}
\right].
\end{gather*}
There is one more simplif\/ication: for an arbitrary $M_{0}=\left[
\begin{matrix}
a_{1} & a_{2}
\\
a_{2} & a_{3}
\end{matrix}
\right] $ if $J=M_{1}^{T}M_{0}M_{1}$ then $ ( J_{11}-J_{22} )
\sin\frac{2\pi\ell}{m}-2J_{12}\cos\frac{2\pi\ell}{m}=2a_{2}$.
So
\begin{gather*}
M_{f}(v) M_{\delta} ( \pi-\psi ) \sim M_{H}\left[
\begin{matrix}
\left( \frac{\psi}{2}\right) ^{\kappa} & -\frac{\delta}{\frac{1}{2}-\kappa }\left( \frac{\psi}{2}\right)^{1+\kappa}
\vspace{1mm}\\
\frac{\delta}{\frac{1}{2}+\kappa}\left( \frac{\psi}{2}\right) ^{1-\kappa} & \left( \frac{\psi}{2}\right)^{-\kappa}
\end{matrix}
\right]  ( M_{1}+\psi\delta M_{2} );
\end{gather*}
and discarding quantities tending to zero (as $\psi\rightarrow0_{+}$) we need to consider the $ ( 1,2 ) $-entry
of
\begin{gather*}
\left[
\begin{matrix}
2^{-\kappa}\psi^{\kappa} & 0
\\
0 & 2^{\kappa}\psi^{-\kappa}
\end{matrix}
\right] M_{H}^{T}\left[
\begin{matrix}
a_{1} & 0
\\
0 & a_{3}
\end{matrix}
\right] M_{H}\left[
\begin{matrix}
2^{-\kappa}\psi^{\kappa} & 0
\\
0 & 2^{\kappa}\psi^{-\kappa}
\end{matrix}
\right],
\end{gather*}
which equals $ ( a_{1}H(\kappa,\delta) -a_{3}H(-\kappa,\delta) )
\frac{\sin\pi\delta}{\cos\pi\kappa}$.
Thus with the normalization constant $c(\kappa,\delta) $
\begin{gather}
K(\phi) =c(\kappa,\delta) L(\phi) ^{T}\left[
\begin{matrix}
H(-\kappa,\delta) & 0
\\
0 & H(\kappa,\delta)
\end{matrix}
\right] L(\phi),
\label{realKmat}
\end{gather}
for $0<\phi=m\theta<\pi$.
The condition that the kernel be positive def\/inite is $H(\kappa,\delta) H(-\kappa,\delta) >0$,
and
\begin{gather*}
H(\kappa,\delta) H(-\kappa,\delta) =\frac{\cos \pi ( \kappa+\delta ) \cos\pi (
\kappa-\delta ) } {\cos^{2}\pi\kappa}=\frac{\sin\pi\left( \frac{\ell}{m}+\kappa\right) \sin \pi\left(
\frac{\ell}{m}-\kappa\right) }{\cos^{2}\pi\kappa},
\end{gather*}
thus $-\frac{\ell}{m}<\kappa<\frac{\ell}{m}$.
In Section~\ref{section6} we show that
\begin{gather*}
c(\kappa,\delta) =\frac{\cos\pi\kappa}{2\pi\cos\pi\delta}.
\end{gather*}

\section{Integrals of harmonic polynomials}  \label{section5}

As in~\cite{Dunkl2013} def\/ine the circle-type pairing by
\begin{gather*}
 \langle f,g \rangle _{S}=\frac{1}{2^{n}n!} \langle f,g \rangle _{G}=\frac{1}{2^{n}n!}\big\langle
e^{\Delta_{\kappa} /2}f,e^{\Delta_{\kappa}/2}g\big\rangle
\end{gather*}
when $\deg f+\deg g=2n$ ($n=0,1,\ldots$); the term ``circle'' refers to the polar
coordinate system.
If~$f$,~$g$ are harmonic then $ \langle f,g \rangle _{G}= \langle f,g \rangle $.
In the coordinate system $ \{ t,\overline{t} \} $ the matrix function (denoted by $K^{\mathbb{C}}$ for this
basis) satisf\/ies $K^{\mathbb{C}}(z) ^{\ast}=K^{\mathbb{C}}(z) $ and $K^{\mathbb{C}}(
zw) =\tau_{\ell}(w) ^{\ast }K^{\mathbb{C}}(z) \tau_{\ell}(w) $ for $w\in W$
(recall $\tau_{\ell}$ is unitary), and $K^{\mathbb{C}}$ is positively homogeneous of degree~$0$.
The Gaussian pairing is
\begin{gather*}
\langle f,g\rangle _{G}=\int_{\mathbb{C}}g(z) K^{\mathbb{C}}(z) f(z)
^{\ast}e^{-\vert z\vert ^{2}/2}dm_{2}(z),
\end{gather*}
where $dm_{2}(z) =dx_{1}dx_{2}$ (for $z=x_{1}+ix_{2}$) and $g(z) $ denotes $ ( g_{1} (
z ),g_{2}(z)  ) $ with $g=g_{1}t+g_{2}\overline{t}$, similarly for~$f$.
In polar coordinates $z=re^{\mathrm{i}\theta}$ the fundamental chamber is $ \big\{  (
r,\theta ):r>0$, $0<\theta<\frac{\pi}{m}\big\} $.
Then
\begin{gather*}
 \langle f,g \rangle _{G}=\sum\limits_{w\in W}\int_{0}^{\infty
}e^{-r^{2}/2}rdr\int_{0}^{\frac{\pi}{m}}g\big( re^{\mathrm{i}\theta }w\big) K^{\mathbb{C}}\big(
e^{\mathrm{i}\theta}w\big) f\big( re^{\mathrm{i}\theta}w\big) ^{\ast}d\theta
\\
\hphantom{\langle f,g \rangle _{G}}{}
=\sum\limits_{w\in W}\int_{0}^{\infty}e^{-r^{2}/2}rdr\int_{0}^{\frac{\pi}{m} }g\big( re^{\mathrm{i}\theta}w\big)
\tau_{\ell}(w) ^{\ast }K^{\mathbb{C}}\big( e^{\mathrm{i}\theta}\big) \tau_{\ell}(w) f\big(
re^{\mathrm{i}\theta}w\big) ^{\ast}d\theta.
\end{gather*}
In particular, if~$f$,~$g$ are homogeneous with $\deg f+\deg g=2n$ then
\begin{gather}
 \langle f,g \rangle _{S}=\frac{1}{2^{n}n!} \langle f,g \rangle _{G}=\sum\limits_{w\in
W}\int_{0}^{\frac{\pi}{m}}g\big( e^{\mathrm{i}\theta}w\big) \tau_{\ell}(w) ^{\ast }K^{\mathbb{C}}\big(
e^{\mathrm{i}\theta}\big) \tau_{\ell}(w) f\big( e^{\mathrm{i}\theta}w\big) ^{\ast}d\theta.
\label{sumgKf}
\end{gather}
Now suppose $f$, $g$ have $m$-parities $r_{1}$, $r_{2}$.
That is, $0\leq r_{1},r_{2}<m$, $f$ is a~sum of monomials $\big\{ z^{a}\overline{z}^{b}
t:a-b+\ell=r_{1}\operatorname{mod}m\big\} $ and $\big\{ z^{a}\overline
{z}^{b}\overline{t}:a-b-\ell=r_{1}\operatorname{mod}m\big\} $ and $g$ has the same property with~$r_{1}$ replaced by~$r_{2}$.
The generic rotation in~$W$ is
\begin{gather*}
\varrho_{j}: \ f ( z,\overline{z},t,\overline{t} ) \mapsto f\big(
z\omega^{j},\overline{z}\omega^{-j},t\omega^{\ell j},\overline{t}\omega^{-\ell j}\big);
\end{gather*}
this maps the monomial $z^{a}\overline{z}^{b}t$ to $\omega^{j ( a-b+\ell ) }z^{a}\overline{z}^{b}t$ and
$z^{a}\overline{z}^{b} \overline{t}$ to $\omega^{j ( a-b-\ell ) }z^{a}\overline{z} ^{b}\overline{t}$.
By hypothesis $f ( e^{\mathrm{i}\theta}\rho _{j} ) \tau_{\ell} ( \rho_{j} ) ^{\ast}=\omega^{r_{1}
j}f ( e^{\mathrm{i}\theta} ) $.
The ref\/lection
\begin{gather*}
\sigma_{j}:\ f ( z,\overline{z},t,\overline{t} ) \mapsto f\big(
\overline{z}\omega^{j},z\omega^{-j},\overline{t}\omega^{\ell j},t\omega^{-\ell j}\big)
\end{gather*}
maps the monomial $z^{a}\overline{z}^{b}t$ to $\omega^{j(a-b+\ell) }z^{b}\overline{z}^{a}\overline{t}=$ and
$z^{a}\overline{z} ^{b}\overline{t}$ to $\omega^{j(a-b-\ell) }z^{b}\overline{z} ^{a}t$.
For a~polynomial $p ( z,\overline{z} ) $ let $p^{\vee }=p ( \overline{z},z ) $.
So if $f=f_{1}t+f_{2}\overline{t}$ has $m$-parity $r_{1}$ then $\sigma_{j}f=\omega^{r_{1}j} ( f_{1}^{\vee
}\overline{t}+f_{2}^{\vee}t ) $.
Break up the sum in formula~\eqref{sumgKf} into two parts, one for rotations (including the identity $\rho_{0}$) and one
for the ref\/lections we obtain
\begin{gather*}
 \langle f,g \rangle _{S}=\sum\limits_{j=0}^{m-1}\omega^{r_{2}j}
\omega^{-r_{1}j}\int_{0}^{\frac{\pi}{m}}g\big( e^{\mathrm{i}\theta}\big) K^{\mathbb{C}}\big(
e^{\mathrm{i}\theta}\big) f\big( e^{\mathrm{i} \theta}\big) ^{\ast}d\theta
\\
\hphantom{\langle f,g \rangle _{S}=}{}
{}+\sum\limits_{j=0}^{m-1}\omega^{r_{2}j}\omega^{-r_{1}j}\int_{0}^{\frac{\pi}{m} }g\symbol{94}\big(
e^{\mathrm{i}\theta}\big) K^{\mathbb{C}}\big( e^{\mathrm{i}\theta}\big) f\symbol{94}\big(
e^{\mathrm{i}\theta}\big) ^{\ast}d\theta,
\end{gather*}
where $g{\symbol{94}}= ( g_{2}^{\vee},g_{1}^{\vee} ) $ and $f{\symbol{94}}= (
f_{2}^{\vee},f_{1}^{\vee} ) $.
The sums vanish unless $r_{1}=r_{2}$ in which case both equal $m$.

Next we convert the matrix in~\eqref{realKmat} to the complex coordinate system.
Recall
\begin{gather*}
K(\phi) =M_{\delta}(\phi) ^{T}M_{f}(v) ^{T}M_{0}M_{f}(v)
M_{\delta}(\phi),
\end{gather*}
in the real coordinates, where
\begin{gather*}
M_{f}(v) =\left[
\begin{matrix}
f_{1}(\kappa,\delta;v) & f_{2}(\kappa,\delta;v)
\\
-f_{2}(-\kappa,\delta;v) & f_{1}(-\kappa,\delta;v)
\end{matrix}
\right],
\qquad
M_{\delta}(\phi) =\left[
\begin{matrix}
\cos\delta\phi & -\sin\delta\phi
\\
\sin\delta\phi & \cos\delta\phi
\end{matrix}
\right].
\end{gather*}
To f\/ind the change-of-basis matrix observe $g_{1}t+g_{2}\overline{t}\mapsto g_{1} ( s_{1}+\mathrm{i}s_{2} )
+g_{2} ( s_{1}-\mathrm{i} s_{2} ) \allowbreak= ( g_{1}+g_{2} ) s_{1}+\mathrm{i} (
g_{1}-g_{2} ) s_{2}$.
Let
\begin{gather*}
B=\left[
\begin{matrix}
1 & \mathrm{i}
\\
1 & -\mathrm{i}
\end{matrix}
\right],
\end{gather*}
then in the complex coordinates
\begin{gather*}
K^{\mathbb{C}}=BM_{\delta}(\phi) ^{T}M_{f}(v) ^{T}M_{0}M_{f}(v) M_{\delta} (
\phi ) B^{\ast}.
\end{gather*}
We f\/ind
\begin{gather*}
BM_{\delta}(\phi) ^{T}=\left[
\begin{matrix}
e^{-\mathrm{i}\delta\phi} & 0
\\
0 & e^{\mathrm{i}\delta\phi}
\end{matrix}
\right] B,
\end{gather*}
and $\delta\phi=\left( \frac{1}{2}-\frac{\ell}{m}\right) m\theta$.
So
\begin{gather*}
K^{\mathbb{C}}=\left[
\begin{matrix}
e^{-\mathrm{i}\delta\phi} & 0
\\
0 & e^{\mathrm{i}\delta\phi}
\end{matrix}
\right] BM_{f}(v) ^{T}M_{0}M_{f}(v) B^{\ast }\left[
\begin{matrix}
e^{\mathrm{i}\delta\phi} & 0
\\
0 & e^{-\mathrm{i}\delta\phi}
\end{matrix}
\right].
\end{gather*}
Then using~\eqref{Mdef}{\samepage
\begin{gather}
K_{11}^{\mathbb{C}}=K_{22}^{\mathbb{C}}=c(\kappa,\delta) G_{1} ( \kappa,\delta,v ),
\nonumber
\\
G_{1} ( \kappa,\delta,v ) :=\big( f_{1}(\kappa,\delta;v) ^{2}+f_{2}(\kappa,\delta;v)
^{2}\big) H(-\kappa,\delta) +\big( f_{1}(-\kappa,\delta;v) ^{2}+f_{2} (
-\kappa,\delta;v ) ^{2}\big) H(\kappa,\delta),
\nonumber
\\
K_{12}^{\mathbb{C}}=\overline{K_{21}^{\mathbb{C}}}=e^{-2\mathrm{i}\delta m\theta}c(\kappa,\delta)
G_{2}(\kappa,\delta;v),\label{G1kdv}
\\
G_{2}(\kappa,\delta;v) :=\big( f_{1}(\kappa,\delta;v) +\mathrm{i} f_{2} (
\kappa,\delta;v ) \big) ^{2}H(-\kappa,\delta) -\big( f_{1}(-\kappa,\delta;v)
+\mathrm{i}~f_{2}(-\kappa,\delta;v) \big) ^{2}H(\kappa,\delta).
\nonumber
\end{gather}
Recall $v=\sin^{2}\frac{\phi}{2}=\frac{1}{2} ( 1-\cos m\theta ) $ and
$e^{-2\mathrm{i}\delta\phi}=e^{-\mathrm{i}\theta ( m-2\ell ) }$ ($=\overline{z}^{m}z^{2\ell}$ on the circle).}

In the special case $f=g=t$ we obtain
\begin{gather}
2=\langle t,t\rangle _{S}=m\int_{0}^{\frac{\pi}{m}}\big( K_{11}^{\mathbb{C}}+K_{22}^{\mathbb{C}}\big)
d\theta=2m\int_{0}^{\frac{\pi }{m}}K_{11}^{\mathbb{C}}d\theta.
\label{normalc}
\end{gather}
In the ``trivial'' case $\kappa=0$ the term $f_{1} ( 0,\delta;v )
^{2}+f_{2} ( 0,\delta;v ) ^{2}=1$ (see Lemma~\ref{f12k0}) and $H(0,\delta) =\cos\pi\delta$ thus
$c(0,\delta) =\frac{1}{2\pi\cos\pi\delta}$.
Consider the evaluation of $\langle f,f\rangle _{S}$ where $f=p_{n}^{(1) }$ and $n=m (
k+1 ) -\ell+r$ with $k\geq0$ and $1\leq r\leq m$; thus
\begin{gather*}
f=z^{m-\ell+r}Q_{k}^{(1) }\left( \kappa,\frac{m-\ell}{m};z^{m},\overline{z}^{m}\right)
t+z^{r}\overline{z}^{m-\ell}Q_{k}^{(2) }\left( \kappa,\frac{m-\ell}{m};z^{m},\overline{z}^{m}\right)
\overline{t}.
\end{gather*}
Because $Q_{k}^{(1) }$ and $Q_{k}^{(2) }$ have real coef\/f\/icients we f\/ind
$f_{1}^{\vee}=\overline{f_{1}}$ and $f_{2}^{\vee }=\overline{f_{2}}$.
The integrand is (suppressing the arguments $\kappa$, $\frac{m-\ell}{m}$ in $Q_{k}^{(j) }$ and omitting the
factor $c(\kappa,\delta) $)
\begin{gather*}
2\big( f_{1}\overline{f_{1}}+f_{2}\overline{f_{2}}\big) G_{1} +2e^{-2\mathrm{i}\delta
m\theta}f_{1}\overline{f_{2}}G_{2}+2e^{2\mathrm{i} \delta m\theta}f_{2}\overline{f_{1}}\overline{G_{2}},
\\
f_{1}\overline{f_{1}}+f_{2}\overline{f_{2}}=Q_{k}^{(1) }\big( z^{m},\overline{z}^{m}\big)
Q_{k}^{(1) }\big( \overline{z}^{m},z^{m}\big) +Q_{k}^{(2) }\big(
z^{m},\overline{z}^{m}\big) Q_{k}^{(2) }\big( \overline {z}^{m},z^{m}\big),
\\
e^{-2\mathrm{i}\delta m\theta}f_{1}\overline{f_{2}}=\overline{z}
^{m}z^{2\ell}z^{m-\ell+r}\overline{z}^{r}z^{m-\ell}Q_{k}^{(1) }\big( z^{m},\overline{z}^{m}\big)
Q_{k}^{(2) }\big( \overline{z}^{m},z^{m}\big)
\\
\hphantom{e^{-2\mathrm{i}\delta m\theta}f_{1}\overline{f_{2}}}{}
=z^{m}Q_{k}^{(1) }\big( z^{m},\overline{z}^{m}\big) Q_{k}^{(2) }\big(
\overline{z}^{m},z^{m}\big),
\\
e^{2\mathrm{i}\delta m\theta}\overline{f_{1}}f_{2}=\overline{z}^{m} Q_{k}^{(1) }\big(
\overline{z}^{m},z^{m}\big) Q_{k}^{(2) }\big( z^{m},\overline{z}^{m}\big).
\end{gather*}
(Note $z\overline{z}=1$ on the circle.) The integrand is now expressed as a~function of~$z^{m}$, that is,
$e^{\mathrm{i}\phi}$.
Similarly suppose $f=p_{n}^{(2) }$ with $n=mk+\ell+r$ so that
\begin{gather*}
p_{n}^{(2) }=z^{r}\overline{z}^{\ell}Q_{k}^{(2) }\left(
\kappa,\frac{\ell}{m};z^{m},\overline{z}^{m}\right) t+z^{\ell+r}Q_{k}^{(1) }\left(
\kappa,\frac{\ell}{m};z^{m},\overline{z}^{m}\right) \overline{t},
\\
e^{-2\mathrm{i}\delta m\theta}f_{1}\overline{f_{2}}=\overline{z}
^{m}z^{2\ell}z^{r}\overline{z}^{\ell}\overline{z}^{\ell+r}Q_{k}^{(2) }\big( z^{m},\overline{z}^{m}\big)
Q_{k}^{(1) }\big( \overline{z}^{m},z^{m}\big)
=\overline{z}^{m}Q_{k}^{(2) }\big( z^{m},\overline{z} ^{m}\big) Q_{k}^{(1) }\big(
\overline{z}^{m},z^{m}\big),
\\
e^{2\mathrm{i}\delta m\theta}\overline{f_{1}}f_{2}=z^{m}Q_{k}^{(2) }\big( \overline{z}^{m},z^{m}\big)
Q_{k}^{(1) }\big( z^{m},\overline{z}^{m}\big).
\end{gather*}

\section{The normalizing constant} \label{section6}

In this section we f\/ind the normalizing constant for $K$ by computing the integral in~\eqref{normalc} and using a~method
that separates the parameters $\delta$ and $\kappa$.
To start we examine the special cases $f_{1}(0,\delta;v) $ and $f_{2}(0,\delta;v)$.

\begin{lemma}
\label{f12k0}
For $-\pi<\phi<\pi$ and any $\delta$
\begin{gather*}
f_{1}\left( 0,\delta;\sin^{2}\frac{\phi}{2}\right)=\cos\delta\phi,
\qquad
f_{2}\left( 0,\delta;\sin^{2}\frac{\phi}{2}\right)=\sin\delta\phi.
\end{gather*}
\end{lemma}

\begin{proof}
The f\/irst part is a~standard formula (see~\cite[(15.4.12)]{Olver2010}).
Brief\/ly, consider $g(v) ={}_{2}F_{1}\left(   \genfrac{}{}{0pt}{}{\delta,-\delta}{1/2} ;v\right) $; it can be
shown directly that $\partial_{\phi}^{2}g(v) =-\delta^{2}g(v)$, $g( 0) =1$ and $g(
1) =\cos\pi\delta$; note $\partial_{\phi}=\sqrt{v(1-v) }\partial_{v}$.
Apply $\partial_{\phi}$ to both sides of the f\/irst formula to obtain
\begin{gather*}
\frac{\delta ( -\delta ) }{1/2}\sqrt{v(1-v) }\;{}_{2}F_{1}\left(
\genfrac{}{}{0pt}{}{1+\delta,1-\delta}{3/2} ;v\right) =-\delta\sin\delta\phi;
\end{gather*}
the left hand side equals $-\delta f_{2}(0,\delta;v) $.
\end{proof}

Next we derive another expression for $K_{11}^{\mathbb{C}}+K_{22}^{\mathbb{C} }$.
Recall $G_{1}(\kappa,\delta;v) $ from~\eqref{G1kdv}.

\begin{lemma}
For $-\frac{1}{2}<\kappa<\frac{1}{2}$ and $\delta>0$
\begin{gather*}
G_{1}(\kappa,\delta;v) =f_{1}(\kappa,\delta;v) f_{1}(-\kappa,\delta;1-v)
+f_{1}(-\kappa,\delta;v) f_{1}(\kappa,\delta;1-v)
\\
\hphantom{G_{1}(\kappa,\delta;v) =}{}
{}-f_{2}(\kappa,\delta;v) f_{2}(-\kappa,\delta;1-v) -f_{2}(-\kappa,\delta;v)
f_{2}(\kappa,\delta;1-v).
\end{gather*}
\end{lemma}

\begin{proof}
From equation~\eqref{fv21mv} we have (replacing $\kappa$ by $-\kappa$ and $v$ by $1-v$)
\begin{gather*}
H(-\kappa,\delta) f_{1}(\kappa,\delta;v)=f_{1}(-\kappa,\delta;1-v)
-\frac{\sin\pi\delta}{\cos\pi\kappa }f_{2}(-\kappa,\delta;v),
\\
H(-\kappa,\delta) f_{1}(\kappa,\delta;v) ^{2}=f_{1}(\kappa,\delta;v) f_{1} (
-\kappa,\delta;1-v ) -\frac{\sin\pi\delta}{\cos\pi\kappa}f_{1}(\kappa,\delta;v) f_{2} (
-\kappa,\delta;v ),
\\
H(\kappa,\delta) f_{1}(-\kappa,\delta;v) ^{2}=f_{1}(-\kappa,\delta;v) f_{1} (
\kappa,\delta;1-v ) -\frac{\sin\pi\delta}{\cos\pi\kappa}f_{1}(-\kappa,\delta;v) f_{2} (
\kappa,\delta;v ),
\end{gather*}
and
\begin{gather*}
H(-\kappa,\delta) f_{2}(\kappa,\delta;v)=\frac{\sin\pi\delta}{\cos\pi\kappa}f_{1} (
-\kappa,\delta;v ) -f_{2}(-\kappa,\delta;1-v),
\\
H(-\kappa,\delta) f_{2}(\kappa,\delta;v) ^{2}=\frac{\sin\pi\delta}{\cos\pi\kappa}f_{1} (
-\kappa,\delta;v ) f_{2}(\kappa,\delta;v) -f_{2}(\kappa,\delta;v) f_{2} (
-\kappa,\delta;1-v ),
\\
H(\kappa,\delta) f_{2}(-\kappa,\delta;v) ^{2}=\frac{\sin\pi\delta}{\cos\pi\kappa}f_{1} (
\kappa,\delta;v ) f_{2}(-\kappa,\delta;v) -f_{2}(-\kappa,\delta;v) f_{2} (
\kappa,\delta;1-v ).
\end{gather*}
Adding up the appropriate equations (the second two in each group) we get the desired result.
Note the cancellations of the terms $f_{1}(\kappa,\delta;v) f_{2}(-\kappa,\delta;v) $ and
$f_{1}(-\kappa,\delta;v) f_{2}(\kappa,\delta;v) $.
\end{proof}

\begin{proposition}
Suppose $-\frac{1}{2}<\kappa<\frac{1}{2}$ and $0\leq\delta<\frac{1}{2}$ then
\begin{gather*}
\int_{0}^{1}G_{1}(\kappa,\delta;v)  ( v(1-v)  )
^{-1/2}dv=\frac{1}{\cos\pi\kappa}\int_{0}^{1}G_{1}(0,\delta;v)  ( v(1-v)  ) ^{-1/2}dv.
\end{gather*}
\end{proposition}

\begin{proof}
By the symmetry $v\rightarrow1-v$ it suf\/f\/ices to integrate $f_{1}(\kappa,\delta;v) f_{1} (
-\kappa,\delta;1-v ) $ and $f_{2}(\kappa,\delta;v) f_{2}(-\kappa,\delta;1-v)$. By def\/inition
\begin{gather*}
f_{1}(\kappa,\delta;v) f_{1}(-\kappa,\delta;1-v) =v^{\kappa}(1-v)
^{-\kappa}\sum\limits_{i,j=0}^{\infty}\frac{(\delta) _{i} ( -\delta ) _{i}(\delta)
_{j} ( -\delta ) _{j}}{i!j!\left( \frac{1}{2}+\kappa\right) _{i}\left( \frac{1}{2}-\kappa\right)
_{j}}v^{i}(1-v) ^{j}.
\end{gather*}
Multiply by $ ( v(1-v)  ) ^{-1/2}$ and integrate over $0<v<1$ to obtain
\begin{gather*}
\sum\limits_{i,j=0}^{\infty}\frac{(\delta) _{i}(-\delta) _{i}(\delta) _{j} (
-\delta ) _{j} }{i!j!\left( \frac{1}{2}+\kappa\right) _{i}\left( \frac{1}{2} -\kappa\right) _{j}}B\left(
i+\frac{1}{2}+\kappa,j+\frac{1}{2}-\kappa\right)
\\
\qquad{} =\Gamma\left( \frac{1}{2}+\kappa\right) \Gamma\left( \frac{1}{2} -\kappa\right)
\sum\limits_{i,j=0}^{\infty}\frac{(\delta) _{i}(-\delta) _{i}(\delta) _{j} (
-\delta ) _{j} }{i!j! ( i+j ) !},
\end{gather*}
because
\begin{gather*}
B\left( i+\frac{1}{2}+\kappa,j+\frac{1}{2}-\kappa\right) =\frac {1}{\Gamma\left( i+j+1\right) }\Gamma\left(
\frac{1}{2}+\kappa\right) \Gamma\left( \frac{1}{2}-\kappa\right) \left( \frac{1}{2}+\kappa\right) _{i}\left(
\frac{1}{2}-\kappa\right) _{j}.
\end{gather*}
The ratio to the same integral with $\kappa$ replaced by $0$ is
\begin{gather*}
\frac{\Gamma\left( \frac{1}{2}+\kappa\right) \Gamma\left( \frac{1} {2}-\kappa\right) }{\Gamma\left( \frac{1}{2}\right)
^{2}}=\frac{1}{\cos \pi\kappa}.
\end{gather*}
Similarly
\begin{gather*}
f_{2}(\kappa,\delta;v) f_{2}(-\kappa,\delta;1-v) =\frac{\delta^{2}}{\left(
\frac{1}{2}+\kappa\right) \left( \frac{1} {2}-\kappa\right) }v^{1+\kappa}(1-v) ^{1-\kappa}
\\
\hphantom{f_{2}(\kappa,\delta;v) f_{2}(-\kappa,\delta;1-v) =}{}
\times\sum\limits_{i,j=0}^{\infty}\frac{(1+\delta) _{i} ( 1-\delta ) _{i} ( 1+\delta )
_{j}(1-\delta) _{j}}{i!j!\left( \frac{3}{2}+\kappa\right) _{i}\left( \frac{3}{2} -\kappa\right)
_{j}}v^{i}(1-v) ^{j}.
\end{gather*}
Multiply by $ ( v(1-v)  ) ^{-1/2}$ and integrate over $0<v<1$ to obtain
\begin{gather*}
\frac{\delta^{2}}{\left( \frac{1}{2}+\kappa\right) \left( \frac{1} {2}-\kappa\right)
}\sum\limits_{i,j=0}^{\infty}\frac{(1+\delta) _{i}(1-\delta) _{i}(1+\delta)
_{j}(1-\delta) _{j}}{i!j!\left( \frac{3}{2}+\kappa\right) _{i}\left( \frac{3}{2}-\kappa\right) _{j}}B\left(
i+\frac{3}{2}+\kappa,j+\frac{3} {2}-\kappa\right)
\\
\qquad{}
=\Gamma\left( \frac{1}{2}+\kappa\right) \Gamma\left( \frac{1}{2} -\kappa\right)
\delta^{2}\sum\limits_{i,j=0}^{\infty}\frac{(1+\delta) _{i}(1-\delta) _{i}(1+\delta)
_{j}(1-\delta) _{j}}{i!j! ( i+j+2 ) !},
\end{gather*}
because
\begin{gather*}
\frac{B\left( i+\frac{3}{2}+\kappa,j+\frac{3}{2}-\kappa\right) }{\left( \frac{1}{2}+\kappa\right) \left(
\frac{1}{2}-\kappa\right) }=\frac{\Gamma\left( \frac{3}{2}+\kappa\right) \Gamma\left( \frac{3} {2}-\kappa\right)
\left( \frac{3}{2}+\kappa\right) _{i}\left( \frac{3} {2}-\kappa\right) _{j}}{\left( \frac{1}{2}+\kappa\right) \left(
\frac {1}{2}-\kappa\right) \Gamma\left( i+j+3\right) }
\\
\hphantom{\frac{B\left( i+\frac{3}{2}+\kappa,j+\frac{3}{2}-\kappa\right) }{\left( \frac{1}{2}+\kappa\right) \left(
\frac{1}{2}-\kappa\right) } }{}
=\Gamma\left( \frac{1}{2}+\kappa\right) \Gamma\left( \frac{1}{2} -\kappa\right) \frac{\left( \frac{3}{2}+\kappa\right)
_{i}\left( \frac {3}{2}-\kappa\right) _{j}}{ ( i+j+2 ) !}.
\end{gather*}
The double series converges by Stirling's formula and the comparison test.
By the same argument as for the f\/irst integral the proof is f\/inished.
\end{proof}

\begin{lemma}
Suppose $0<\delta<\frac{1}{2}$ then
\begin{gather*}
\int_{0}^{1}f_{1}(0,\delta;v) f_{1} ( 0,\delta;1-v ) \frac{dv}{\sqrt{v(1-v) }}=\frac{\pi}{2}\cos\delta\pi +\frac{1}{2\delta}\sin\pi\delta,
\\
\int_{0}^{1}f_{2}(0,\delta;v) f_{2} ( 0,\delta;1-v ) \frac{dv}{\sqrt{v(1-v) }}
=-\frac{\pi}{2}\cos\delta\pi +\frac{1}{2\delta}\sin\pi\delta.
\end{gather*}

\end{lemma}

\begin{proof}
For $v=\sin^{2}\frac{\phi}{2}$ one has $dv=\sqrt{v(1-v) }d\phi$.
Also $1-v=\cos^{2}\frac{\phi}{2}=\sin^{2}\frac{\pi-\phi}{2}$, thus by Lemma~\ref{f12k0} the two integrals equal
\begin{gather*}
\int_{0}^{\pi}\cos\delta\phi~\cos\delta ( \pi-\phi ) d\phi=\frac{\pi}{2}\cos\delta\pi+\frac{1}{2\delta}\sin\delta\pi,
\\
\int_{0}^{\pi}\sin\delta\phi~\sin\delta ( \pi-\phi ) d\phi=-\frac{\pi}{2}\cos\delta\pi+\frac{1}{2\delta}\sin\delta\pi,
\end{gather*}
respectively.
\end{proof}

Combining the lemma and the proposition we prove:

\begin{proposition}
Suppose $-\frac{1}{2}<\kappa<\frac{1}{2}$ and $0\leq\delta<\frac{1}{2}$ then
\begin{gather*}
\int_{0}^{1}G_{1}(\kappa,\delta;v)  ( v(1-v) )
^{-1/2}dv=2\pi\frac{\cos\delta\pi}{\cos\kappa\pi}.
\end{gather*}

\end{proposition}

To determine the normalizing constant: recall the requirement from~\eqref{normalc}
\begin{gather*}
2=\langle t,t\rangle _{S}=2m\int_{0}^{\frac{\pi}{m}}
K_{11}^{\mathbb{C}}d\theta=2\int_{0}^{\pi}K_{11}^{\mathbb{C}}d\phi
\\
\hphantom{2}{}
=2c(\kappa,\delta) \int_{0}^{1}G_{1}(\kappa,\delta;v)  ( v(1-v)  )
^{-1/2}dv=4\pi c(\kappa,\delta) \frac{\cos\delta\pi}{\cos\kappa\pi},
\\
c(\kappa,\delta)=\frac{1}{2\pi}\frac{\cos\kappa\pi} {\cos\delta\pi}.
\end{gather*}
This implies
\begin{gather*}
\det K=\frac{\sin\pi\left( \frac{\ell}{m}+\kappa\right) \sin\pi\left( \frac{\ell}{m}-\kappa\right)
}{4\pi^{2}\sin^{2}\frac{\pi\ell }{m}}=\frac{1}{4\pi^{2}}\left( 1-\left( \frac{\sin\pi\kappa}{\sin\frac
{\pi\ell}{m}}\right) ^{2}\right)
\end{gather*} (and $\det K^{\mathbb{C}}= \vert \det B \vert ^{2}\det K=4\det K$).
The weight function $K$ is actually integrable for $-\frac{1}{2}<\kappa<\frac{1}{2}$ but does not give a~positive form
when $\frac{\ell}{m}\leq \vert \kappa \vert <\frac{1}{2}$.

\pdfbookmark[1]{References}{ref}
\LastPageEnding

\end{document}